\documentclass[twocolumn]{autart}

\usepackage{amsmath}
\usepackage{amssymb}
\usepackage{color}
\usepackage{graphicx}
\graphicspath{ {pics/} }
\usepackage{mathptmx}
\usepackage{subcaption}
\usepackage{wrapfig}
\usepackage{units}
\allowdisplaybreaks

\begin{document}
\begin{frontmatter}
	\title{Output Feedback Control of Two-lane Traffic Congestion}

\thanks[*]{Corresponding author.}
	
\author[*]{Huan Yu}\ead{huy015@eng.ucsd.edu}\;\text{and}\; 
\author{Miroslav Krstic}\ead{krstic@ucsd.edu}               

\address{Department of Mechanical and Aerospace Engineering, University of California, San Diego, La Jolla, CA, 92093, United States}  

                        
\begin{abstract}
	This paper develops output feedback boundary control to mitigate traffic congestion of a unidirectional two-lane freeway segment. The macroscopic traffic dynamics are described by the Aw-Rascle-Zhang (ARZ) model respectively for both the fast and slow lanes. The traffic density and velocity of each of the two lanes are governed by coupled $2 \times 2$ nonlinear hyperbolic partial differential equations (PDEs). Lane-changing interactions between the two lanes lead to exchanging source terms between the two pairs second-order PDEs. Therefore, we are dealing with $4 \times 4$ nonlinear coupled hyperbolic PDEs. Based on driver's preference for the slow and fast lanes, a reference system of lane-specific uniform steady states in congested traffic is chosen. To stabilize traffic densities and velocities of both lanes to the steady states, two distinct variable speed limits (VSL) are applied at outlet boundary, controlling the traffic velocity of each lane. Using backstepping transformation, we map the coupled heterodirectional hyperbolic PDE system into a cascade target system, in which traffic oscillations are damped out through actuation of the velocities at the downstream boundary. Two full-state feedback boundary control laws are developed. We also design a collocated boundary observer for state estimation with sensing of densities at the outlet. Output feedback boundary controllers are obtained by combining the collocated observer and full-state feedback controllers. The finite time convergence to equilibrium is achieved for both the controllers and observer designs. Numerical simulations validate our design in two different traffic scenarios.
\end{abstract}

\end{frontmatter}

\section{Introduction}

Traffic congestion on freeways been investigated intensively over the past decades. Motivations behind are to understand the formation of traffic congestion, and further to prevent or suppress instabilities of traffic flow. Macroscopic modeling of traffic dynamics is to describe evolution of aggregated traffic state values including traffic density and velocity. Traffic dynamics are governed by hyperbolic PDEs, including the first-order model by Ligthill, Whitham and Richards (LWR), the second-order Payne-Whitham model, and the second-order Aw-Rascle-Zhang (ARZ) model~\cite{AW:00}~\cite{Zhang}. The LWR model is a conservation law of traffic density. It is simple yet powerful in understanding the formation and propagation of traffic shockwaves on freeway. But it fails to describe stop-and-go traffic, also known as "jamiton"~\cite{Seibold:09}. The oscillations of densities and velocities travel with traffic stream, causing unsafe driving conditions, increased consumptions of fuel and delay of travel time.

In order to describe this common phenomenon, the second-order traffic models are developed, consisting of nonlinear hyperbolic PDEs of traffic density and velocity. The ARZ model advances the PW model in correctly prediction of propagation of traffic velocity. Validation of the ARZ model with empirical traffic data is conducted in~\cite{Fan:13}. ~\cite{Belletti:15} discusses the heterodirectional propagations of characteristics waves for congested traffic by the ARZ model. Previous work of authors~\cite{Yu:19}~\cite{Yu:18} discuss the linear stability of uniform steady states of the nonlinear hyperbolic ARZ model. Instabilities appear in the congested regime of the ARZ model when drivers on the road are aggressive. To tackle the traffic congestion, traffic management infrastructures like VSL and ramp metering provide opportunities for boundary actuation of traffic velocity and flow. Many recent efforts~\cite{Bekiaris-Liberis:18}~\cite{Karafyllis:18}~\cite{Yu:19}~\cite{Yu:18}~\cite{Liguo:19}~\cite{Liguo:18}~\cite{Liguo:17} are focused on boundary control of traffic PDE model. Among these results, \cite{Yu:19} firstly applied the backstepping control design to stabilize the stop-and-go traffic of the one-lane ARZ model with ramp metering.

The above-mentioned models treat multi-lane freeway traffic cumulatively as a single lane by assuming averaged velocity and density over cross section of all lanes. The individual dynamics of each lane and inter-lane interactions are neglected. However, distinct density and velocity equilibrium exist in multi-lane problems. The differences of velocities give rise to lane-changing interactions and further lead to traffic congestion~\cite{Daganzo}. To address the phenomenon, a number of macroscopic multi-lane models~\cite{Dirk:97}~\cite{Herty:03}~\cite{PG}~\cite{Klar1:98}~\cite{Klar2:03} have been developed from microscopic, then kinetic to macroscopic descriptions. In this paper, we adopt the multi-lane ARZ traffic model proposed by~\cite{Klar1:98}~\cite{Herty:03} to describe a two-lane freeway traffic with lane-changing between the two lanes. Lane interactions appear as interchanging source terms in the system, leading to more involved couplings and a higher order of PDEs. The complexity of the multi-lane model is greatly increased compared to the one-lane problem.

Feedback boundary control design for a general class of hyperbolic PDEs using backstepping method are studied in ~\cite{Anfinsen:17}~\cite{Jean:16}~\cite{Coron:13}~\cite{DiMeglio:13}~\cite{Deutscher:2017}~\cite{long}~\cite{Yu:19}. In~\cite{DiMeglio:13}, stabilization of a $n+1$ counter convecting hyperbolic PDEs is achieved with a single boundary. ~\cite{long} presents a solution to output feedback of a fully general case of heterodirectional $n+m$ first-order linear coupled hyperbolic PDEs. Actuation of all the $m$ PDEs from the same boundary is required to stabilize the system in finite time. A shorter convergence time is further obtained in \cite{Jean:16} by modifying the target system structure.

In this problem, a two-lane ARZ model of a freeway segment presents heterodirectional $2+2$ coupled nonlinear hyperbolic PDEs, governing the traffic densities and velocities of the fast and slow lanes. We aim to stabilize the oscillations in the two-lane traffic using the PDE backstepping method, based on the stabilization results in~\cite{long}. Actuation of traffic velocities at the outlet boundary are realized by two VSLs.

The contribution of this paper is the following. This is the first result dealing with control problem of multi-lane traffic PDE model. The dynamics of two-lane traffic are studied from control perspectives. Theoretical result of output feedback control of the general class of heterodirectional linear hyperbolic PDE systems is developed in~\cite{Jean:16}~\cite{long}, but has never been applied in traffic application. Being the first paper to adopt the methodology, our result opens the door for solving related multi-lane traffic problems with PDE control techniques. Furthermore, we advance the theoretical results in~\cite{Jean:16}~\cite{long} by proposing a collocated boundary observer and controller design. The output feedback controllers in both papers are constructed with full-state feedback controllers and an anti-collocated observer. In implementation, collocated boundary observer and controllers are more practically applicable. We bridge this gap by developing a observer with sensing at outlet, which is also a more challenging problem in the design for the system of this paper.

The paper is organized as follows: in Section 2 we introduce the two-lane ARZ traffic model and then linearize the nonlinear hyperbolic PDEs around uniform steady states. In Section 3 backstepping transformation is derived for the linearized model in Riemann coordinates. We present full-state feedback control laws to actuate outlet boundary velocities. In Section 4, we design collocated boundary observers and then obtain output feedback control laws. In Section 5, control design in two different traffic scenarios are discussed and tested with numerical simulation. 

\section{Problem Statement}	
In this section, two-lane traffic ARZ model is introduced. We derive lane-specific uniform steady states according to the drivers' overall preference for the lanes and then linearize the nonlinear system around the steady states. The linearized system is transformed to Riemann variables. A coupled $4 \times 4$ first-order hetero-directional hyperbolic system is obtained for control design.

\subsection{Two-lane traffic ARZ model}
	\begin{figure}[t!]
		\centering					
		\includegraphics[width=2.5cm]{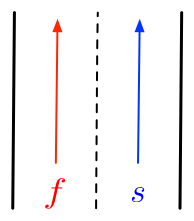}	
		\caption{A unidirectional freeway segment of the fast and slow lanes }
	\end{figure}
The two-lane traffic on unidirectional roads is described with the following two-lane traffic ARZ model by ~\cite{Klar1:98}~\cite{Herty:03}. The diagram in Fig.1 is shown with the faster lane on the left and slower lane on the right. The two-lane traffic ARZ model is given by
\begin{align}
		\partial_t \rho_{ f}+ \partial_x( \rho_{ f} v_{ f})=& \frac{1}{T_{ s}}\rho_ { s} - \frac{1}{T_{ f}}\rho_{ f}, \label{rv1}\\
\notag		\partial_t ( \rho_{ f} v_{ f})+\partial_x( \rho_{ f} v_{ f}^2)-(\gamma p_{ f})\partial_x v_{ f}=& \frac{1}{T_{ s}}\rho_ { s}v_ { s} - \frac{1}{T_{ f}}\rho_{ f}v_ { f}\\& + \frac{\rho_ { f}(V(\rho_ { f})-v_{ f})}{T_{ f}^e},\\
			\partial_t \rho_{ s}+ \partial_x( \rho_{ s} v_{ s})=&\frac{1}{T_{ f}}\rho_ { f} - \frac{1}{T_{ s}}\rho_{ s}, \\
\notag	\partial_t ( \rho_{ s} v_{ s})+\partial_x( \rho_{ s} v_{ s}^2)-(\gamma p_{ s})\partial_x v_{ s}=& \frac{1}{T_{ f}}\rho_ { f}v_ { f} -\frac{1}{T_{ s}}\rho_{ s}v_ { s} \\&+ \frac{\rho_ { s}(V(\rho_ { s})-v_{ s})}{T_{ s}^e}.\label{rv4}
	\end{align}	
The traffic density $\rho_i(x,t)$ and velocity $v_i(x,t)$ $ (i=  f, s)$ are defined in $x\in[0,L]$, $t\in[0,\infty)$, where $L$ is the length of the freeway segment. 
The above nonlinear hyperbolic PDEs consist of two subsystems of second-order nonlinear hyperbolic PDEs, each describing one-lane traffic dynamics. Lane-changing interactions and drivers' behavior adapting to the traffic appear as source terms on the right hand side of PDEs.    

The variable $p_i(\rho_i)$ is defined as the traffic density pressure
\begin{align}
p_i(\rho_i)=v_{m}\left(\frac{\rho_i}{\rho_m}\right)^\gamma,
\end{align}
which is an increasing function of density $\rho_{ i}$.  $v_{m}$ is the maximum traffic velocity, $\rho_m$ is the maximum traffic density and the constant coefficient $\gamma \in \mathbb{R}_{+}$ reflects the aggressiveness of drivers on road. The parameter $T^e_i$ is defined as relaxation time that reflects driver's behavior adapting to the traffic equilibrium velocity in the lane $i$. The parameter $T_{i}$ describes the driver's preference for remaining in lane $i$, which relates to the both lanes' density and velocity. We consider them to be constant coefficients in this paper.

The equilibrium velocity-density relationship $V(\rho)$ is given in the form of the Greenshield's model,
	\begin{align}
		V(\rho_i)=v_m\left(1-\left(\frac{\rho_i}{\rho_m}\right)^{\gamma}\right). \label{green}
	\end{align}
We choose the Greenshield's model for $V(\rho)$ due to its simplicity but the control design presented later is not limited by this choice. Note that the equilibrium velocity-density model \eqref{green} is for cumulative single lane traffic.  Distinct velocity equilibrium does exist in each of the two lanes~\cite{Klar2:03}. The lane-specific steady traffic velocities will be discussed in the following section.
	\subsection{Driver's preference for two lanes}
We consider to linearize the nonlinear hyperbolic system $(\rho_i, v_i)$ around uniform steady states $(\rho_{i}^\star, v_{i}^\star)$. We obtain the following equations
	\begin{align}
	    \frac{1}{T_{ s}}\rho_ { s}^\star - \frac{1}{T_{ f}}\rho_{ f}^\star=&0, \label{d1}\\
		\frac{1}{T_{ s}}\rho_ { s}^\star v_ { s}^\star - \frac{1}{T_{ f}}\rho_{ f}^\star v_ { f}^\star + \frac{\rho_ { f}^\star(V(\rho_ { f}^\star)-v_{ f}^\star)}{T_{ f}^e} =& 0,\\
		\frac{1}{T_{ f}}\rho_ { f}^\star v_ { f}^\star -\frac{1}{T_{ s}}\rho_{ s}^\star v_ { s}^\star + \frac{\rho_ { s}^\star (V(\rho_ { s}^\star )-v_{ s}^\star )}{T_{ s}^e} =& 0.
	\end{align}
	The steady state density-velocity relations are defined based on \eqref{green}.  Thus the steady states $(\rho_{ f}^\star, v_{ f}^\star, \rho_{ s}^\star, v_{ s}^\star)$ need to satisfy
	\begin{align}
		\rho_{ f}^\star =& \sigma \rho_{ s}^\star , \label{s1}\\
		v_{ f}^\star =& v_m \left(1-r_f\left(\frac{\rho_f^\star}{\rho_m}\right)^{\gamma}\right),  \\
		v_{ s}^\star =& v_m \left(1-r_s\left(\frac{\rho_s^\star}{\rho_m}\right)^{\gamma}\right), \label{s3}
	\end{align}
	where $v_{ f}^\star$ and $v_{ s}^\star$ differ from single-lane $V(\rho_i)$. The ratio coefficients $r_f$ and $r_s$ are defined as 
	\begin{align}
		r_f = & \frac{1+\left(\frac{1}{\sigma}\right)^{\gamma} \frac{T_f^e}{T_f}+\frac{T_s^e}{T_s}}{1+\frac{T_f^e}{T_f}+\frac{T_s^e}{T_s}},\\
		r_s = & \frac{1+ \frac{T_f^e}{T_f}+\frac{T_s^e}{T_s}\left({\sigma}\right)^\gamma}{1+\frac{T_f^e}{T_f}+\frac{T_s^e}{T_s}}.
	\end{align}	
The parameter $\sigma$ defines driver's preference for the fast lane over slow lane according to \eqref{d1},\eqref{s1},
	\begin{align}
		\sigma = \frac{T_f}{T_s}.
	\end{align}
Compared with the single-lane Greenshield's model in \eqref{green}, the relations of steady state traffic velocities $v_{ i}^\star$ and densities $\rho_{ i}^\star$ depend on the drivers lane-changing preference parameter $\sigma$.

\begin{figure}[t!]	
	\begin{subfigure}{0.24\textwidth}
		\includegraphics[width=1\linewidth]{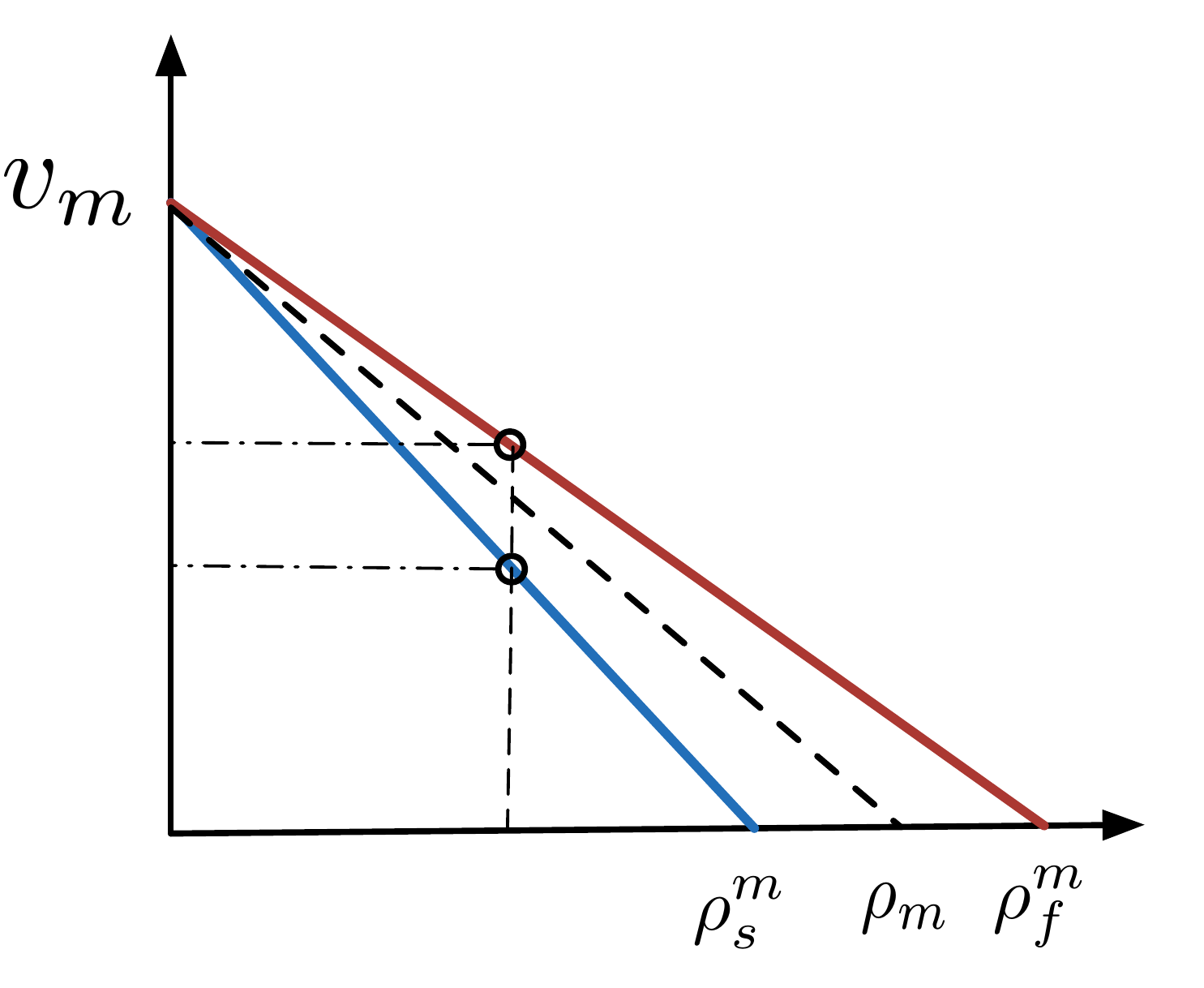} 
		\caption{velocity-density relation}
		\label{fig:subim1}
	\end{subfigure}
	\begin{subfigure}{0.23\textwidth}
		\includegraphics[width=1\linewidth]{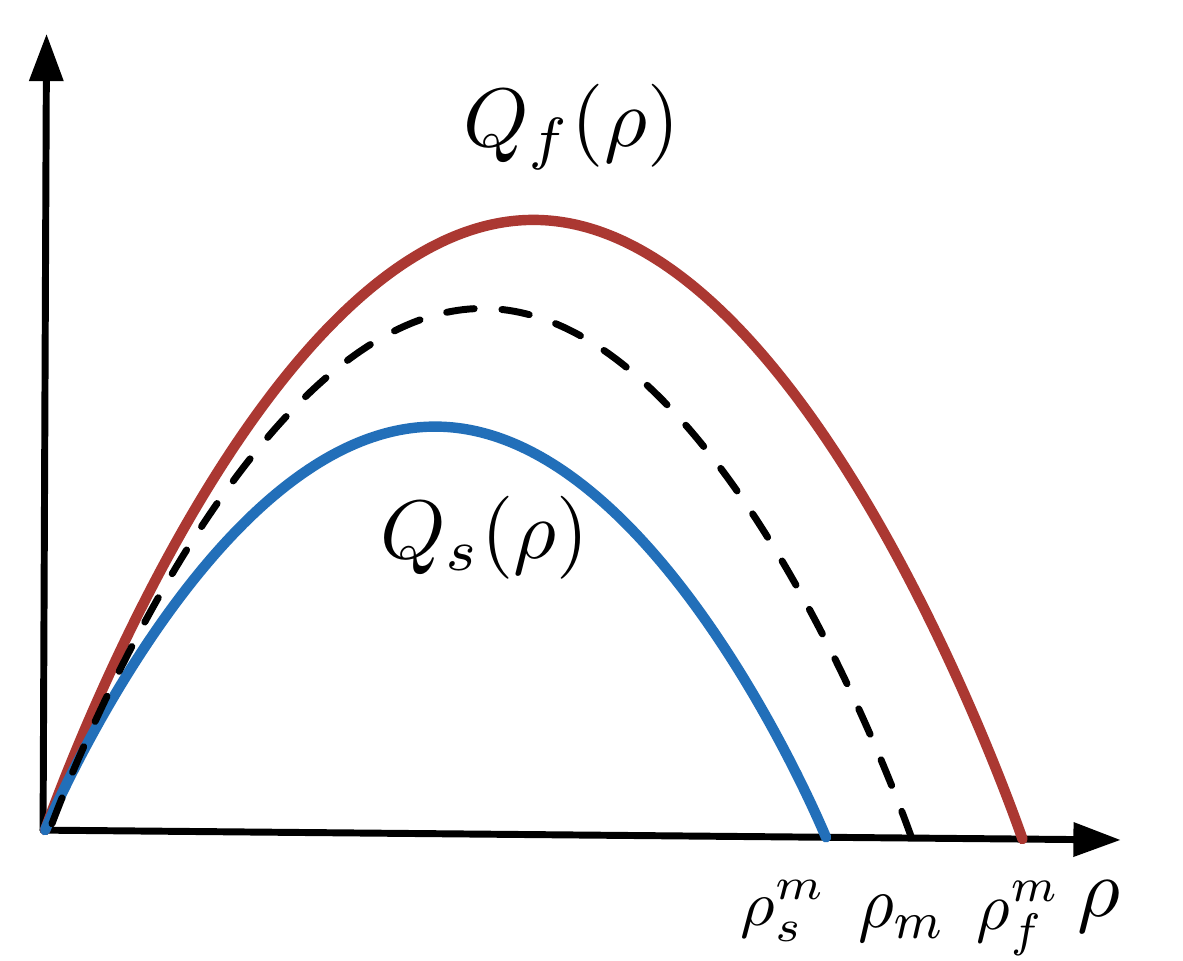}
		\caption{fundamental diagram}
		\label{fig:subim2}
	\end{subfigure}	
	\caption{Steady states of one lane, fast and slow lane}
	\label{fig:image2}
\end{figure}
Assuming that overall drivers prefer fast lane over slow lane, we use Fig.2 $(\sigma >1, \gamma=1)$ to show the equilibrium velocity-density relation and fundamental diagram of the single-lane, the fast and slow lane. $\rho_m$ represents the equivalent maximum density of the single-lane. The actual maximum density in the fast lane $\rho_{f}^m$ and in the fast lane $\rho_{s}^m$ are related to $\rho_m$ by
\begin{align}
	\rho_{f}^m =& \sqrt[\gamma]{r_f} \rho_m,\\
	\rho_{s}^m =& \sqrt[\gamma]{r_s} \rho_m.
\end{align}

$\bullet \quad   \sigma >1 \implies r_f< 1 < r_s $
	
	 If drivers prefer the fast lane, the decrease of velocity gets steeper in the slow lane and less steep in the fast lane. At the same density, the fast lane traffic is "more tolerant to risk" of high density than in the single-lane case, and the slow lane traffic is "less tolerant to risk" than in the single-lane case. As a result, the traffic flux of fast lane is higher than the slow lane at the same density in the fundamental diagram shown in fig.2.

$\bullet \quad   \sigma <1 \implies r_f > 1 >r_s $	

Drivers prefer the slow lane. The decrease of velocity is steeper in the fast lane than that of slow lane at the same density. The slow lane is more tolerant to high density and the traffic flux is higher in the slow lane.

In general, the activities of lane changing segregate the drivers into the more "risk-tolerant" ones in the fast lane and the more "risk-averse" in the slow lane. The risk-tolerant drivers prefer to drive with a faster speed at the same density, compared with risk-averse drivers.

\subsection{Linearized two-lane ARZ model}
Before linearizing the nonlinear system \eqref{rv1}-\eqref{rv4} to steady states \eqref{s1}-\eqref{s3}, we consider the following boundary conditions of  $(\rho_i, v_i)$-system. We assume constant traffic flux entering from the inlet boundary $x=0$ of the two lanes.
\begin{align}
	 q_i^\star=\rho^\star_i v_i^\star.
\end{align} 
Two VSLs implemented at the outlet $U_f(t)$ and $U_s(t)$ actuate the traffic velocity variations for the fast and slow lanes respectively.
\begin{align}
\rho_{ f }(0,t)=&\frac{\rho^\star_f v_f^\star}{v_f(0,t)},\label{bc1}\\
v_{ f }(L,t) =& U_f(t)+v_f^\star,\\
\rho_{ s }(0,t)=&\frac{\rho^\star_s v_s^\star}{v_s(0,t)},\\
v_{ s}(L,t) =& U_s(t)+v_s^\star,\label{bc4}
\end{align}	
Then we linearize the above nonlinear hyperbolic system $(\rho_{ f}, v_{ f}, \rho_{ s}, v_{ s})$ around steady states $(\rho_{ f}^\star, v_{ f}^\star, \rho_{ s}^\star, v_{ s}^\star)$ that satisfy \eqref{s1}-\eqref{s3}.
The deviations from the steady states are defined as
\begin{align}
\tilde \rho_{ f} = &\rho_{ f} - \rho_{ f}^\star, \quad \tilde v_{ f} = v_{ f} - v_{ f}^\star,\\
\tilde \rho_{ s} =& \rho_{ s} - \rho_{ s}^\star, \quad \tilde v_{ s} = v_{ s} - v_{ s}^\star.		
\end{align}
The linearized hyperbolic system is obtained 
\begin{align}
\partial_t  \tilde \rho_{ s}  +  v_s^\star 	\partial_x  \tilde \rho_{ s} + \rho_s^\star \partial_x  \tilde v_{ s}=&
 - \frac{1}{T_{ s}} \tilde \rho_{ s} + 		\frac{1}{T_{ f}} \tilde \rho_ { f}, \\	
\partial_t  \tilde \rho_{ f}  +  v_f^\star 	\partial_x  \tilde \rho_{ f} + \rho_f^\star \partial_x  \tilde v_{ f}=&
\frac{1}{T_{ s}} \tilde \rho_ { s} - \frac{1}{T_{ f}} \tilde \rho_{ f}, \\
\notag		\partial_t \tilde v_{ s} + (v_s^\star - \gamma p_s^\star) \partial_x  \tilde v_{ s} 
=& -\frac{1}{T_{ s}} \frac{v_{ f}^\star - v_{ s}^\star }{\rho_ s^\star}\tilde \rho_{ s}+ \frac{1}{T_{ f}} \frac{v_{ f}^\star - v_{ s}^\star}{\rho_ { s}^\star}\tilde \rho_{ f }\\ & +  \frac{1}{T_{s}} (\tilde v_{ f}- \tilde v_{ s} ) + \frac{\tilde \rho_{ s} V'(\rho_ { s}^\star)- \tilde v_{ s}}{T_{ s}^e},\\	
\notag		\partial_t \tilde v_{ f} + (v_f^\star - \gamma p_f^\star) \partial_x  \tilde v_{ f} =
& \frac{1}{T_{ s}} \frac{v_{ s}^\star - v_{ f}^\star}{\rho_ { f}^\star} \tilde \rho_{ s }-\frac{1}{T_{ f}} \frac{v_{ s}^\star - v_{ f}^\star}{\rho_ { f}^\star} \tilde \rho_{ f}\\ &+  \frac{1}{T_{ f}} (\tilde v_{ s}- \tilde v_{ f} )  + \frac{\tilde \rho_{ f} V'(\rho_ { f}^\star)- \tilde v_{ f}}{T_{ f}^e},
\end{align}
with the linearized boundary conditions
\begin{align}
	\tilde \rho_s(0,t) =& -\frac{\rho_{ s}^\star}{v_{ s}^\star}\tilde v_s(0,t),\\
	\tilde \rho_f(0,t) =& -\frac{\rho_{ f}^\star}{v_{ f}^\star}\tilde v_f(0,t),\\
	\tilde v_{ s }(L,t) =& U_s(t),\\
	\tilde v_{ f }(L,t) =& U_f(t).
\end{align}
In order to diagonalize the spatial derivatives on the left hand side of the equations, we write the above linearized hyperbolic system in the Riemann coordinates $(\tilde w_f, \tilde v_f, \tilde w_s, \tilde v_s  )$ as
\begin{align}
	\tilde w_s =& \frac{\gamma p_s^\star}{\rho_s^\star} \tilde \rho_{ s}  + \tilde v_{ s}, \quad \tilde v_s = \tilde v_s, \label{rt2}\\
    \tilde w_f =& \frac{\gamma p_f^\star}{\rho_f^\star} \tilde \rho_{ f}  + \tilde v_{ f}, \quad \tilde v_f = \tilde v_f.\label{rt1}
\end{align}
We consider the congested regime in~\cite{Yu:19} where steady state traffic density disturbances convect downstream and the velocity disturbances travel upstream. Therefore the following conditions hold for the characteristic speeds of $\tilde v_i$,
\begin{align}
v_s^\star - \gamma p_s^\star<0, \quad v_f^\star - \gamma p_f^\star<0.  
\end{align}
We obtain a coupled $4 \times 4$ first-order hetero-directional hyperbolic system in $( \tilde w_s, \tilde w_f,\tilde v_s, \tilde v_f)$,
\begin{align}
	\partial_t \tilde w_s + v_s^\star \partial_x \tilde w_s  =& a^{ww}_{11} \tilde w_s + a^{ww}_{12}\tilde w_f + a^{wv}_{11} \tilde v_{ s} + a^{wv}_{12} \tilde v_{ f},\label{wv1}\\
	\partial_t \tilde w_f + v_f^\star \partial_x \tilde w_f  =& a^{ww}_{21} \tilde w_s + a^{ww}_{22}\tilde w_f + a^{wv}_{21} \tilde v_{ s} + a^{wv}_{22} \tilde v_{ f}, \\
	\partial_t \tilde v_{ s} - (\gamma p_s^\star - v_s^\star) \partial_x  \tilde v_{ s} =& a^{vw}_{11} \tilde w_s + a^{vw}_{12}\tilde w_f + a^{vv}_{11} \tilde v_{ s} + a^{vv}_{12} \tilde v_{f},\\
	\partial_t \tilde v_{ f} - (\gamma p_f^\star-v_f^\star ) \partial_x  \tilde v_{ f} =& a^{vw}_{21}\tilde w_s + a^{vw}_{22} \tilde w_{ f} + a^{vv}_{21} \tilde v_s  + a^{vv}_{22} \tilde v_{f},\\
		\tilde w_s(0,t) =&  k_s \tilde v_s(0,t),\\
		\tilde w_f(0,t) =& k_f \tilde v_f(0,t),\\
		\tilde v_{s}(L,t) =& U_s(t),\\
		\tilde v_{f}(L,t) =& U_f(t), \label{wv2}
\end{align}
where the constant boundary coefficients $k_i$ are defined as
\begin{align}
	k_i =  -\frac{\gamma p_i^\star-v_{ i}^\star}{v_{ i}^\star},
\end{align}
and the constant parameter block matrix $\{A\}$ is denoted by
\begin{align}
	A = \begin{bmatrix}
	A^{ww} & A^{wv}\\
    A^{vw} & A^{vv}\\	
	\end{bmatrix}.
\end{align}
The elements of sub-matrices of $\{A\}$ are defined as, 
\begin{align}
    a^{ww}_{11} = &  - \frac{1}{T_s^e} - \frac{1}{T_s}\frac{v_{ f}^\star - v_{ s}^\star + \gamma p_{ s }^\star}{\gamma p_{ s }^\star},\quad a^{ww}_{12} = \frac{1}{T_s}\frac{v_f^\star-v_s^\star+\gamma p_{s}^\star}{\gamma p_{f}^\star},\\
     a^{ww}_{21} = & \frac{1}{T_f}\frac{v_s^\star-v_f^\star+\gamma p_{ f }^\star}{\gamma p_{ s }^\star}, \quad a^{ww}_{22} =  - \frac{1}{T_f^e} - \frac{1}{T_f}\frac{v_{ s}^\star - v_{ f}^\star + \gamma p_{ f }^\star}{\gamma p_{ f }^\star}, \\
    a^{wv}_{11} = &  \frac{1}{T_s}\frac{v_{ f}^\star - v_{ s}^\star}{\gamma p_{s}^\star}, \quad a^{wv}_{12} = -\frac{1}{T_s}\frac{(\gamma p_{s}^\star-v_s^\star) - (\gamma p_{ f }^\star-v_f^\star)}{\gamma p_{f}^\star} , \\ 
    a^{wv}_{21}= & -\frac{1}{T_f}\frac{(\gamma p_{ f }^\star-v_f^\star)-(\gamma p_{s}^\star-v_s^\star)}{\gamma p_{ s }^\star}, \quad a^{wv}_{22} =   \frac{1}{T_f}\frac{v_{ s}^\star - v_{ f}^\star }{\gamma p_{ f }^\star}, \\
    a^{vw}_{11} = & - \frac{1}{T_s^e}  - \frac{1}{T_s}\frac{v_{ f}^\star - v_{s}^\star}{\gamma p_{s}^\star},\quad a^{vw}_{12} =  \frac{1}{T_s}\frac{v_f^\star-v_s^\star}{\gamma p_{f}^\star},\\
    a^{vw}_{21} = & \frac{1}{T_f}\frac{v_s^\star-v_f^\star}{\gamma p_{ s }^\star}, \quad a^{vw}_{22} =  - \frac{1}{T_f^e}  - \frac{1}{T_f}\frac{v_{ s}^\star - v_{ f}^\star}{\gamma p_{ f }^\star},\\
    a^{vv}_{11} = &\frac{1}{T_s}\frac{v_f^\star-v_s^\star-\gamma p_{s}^\star}{\gamma p_{s}^\star}, \quad a^{vv}_{12} =  - \frac{1}{T_s}\frac{v_{f}^\star - v_{s}^\star - \gamma p_{f}^\star}{\gamma p_{f}^\star},\\
    a^{vv}_{21} = &-\frac{1}{T_f}\frac{v_s^\star-v_f^\star-\gamma p_{ s }^\star}{\gamma p_{ s }^\star},\quad a^{vv}_{22} =   - \frac{1}{T_f}\frac{v_{ s}^\star - v_{ f}^\star - \gamma p_{ f }^\star}{\gamma p_{ f }^\star}.
\end{align}
The flow diagram of $(\tilde w_i, \tilde v_i)$-system is shown in Fig.3. The $4 \times 4$ first-order hyperbolic system is composed of two coupled second-order heterodirectional hyperbolic systems. States $\tilde w_i$ convect downstream while states $\tilde v_i$ propagate upstream. We use two VSLs to damp out the oscillations to zero from the outlet. 
\begin{figure}[t!]
	\centering					
	\includegraphics[width=6cm]{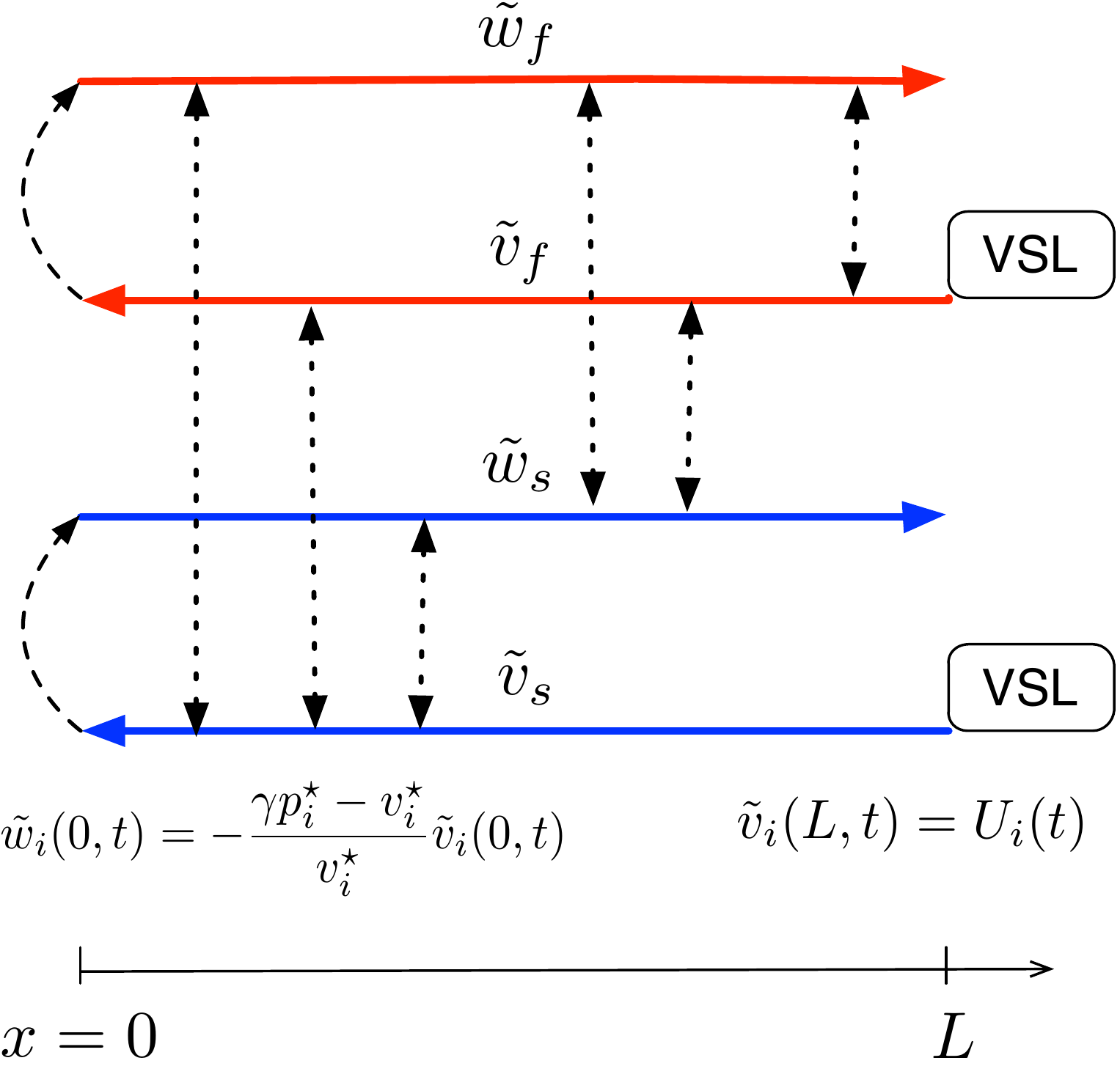}	
	\caption{Flow diagram of linearized two-lane ARZ model}
\end{figure}

\section{Full-state Feedback Control Design with VSLs}	
To apply the backstepping approach and to design boundary control for the system in \eqref{wv1}-\eqref{wv2}, we scale the state variables $\tilde v_s$ and $\tilde v_f$ in space to cancel the diagonal terms in their equations. The Riemann variables $\tilde w_s$ and $\tilde w_f$ remain to be the same. The scaled variables $\bar v_s$ and $\bar v_f$ are defined as
\begin{align}
\quad  \bar v_s =& \exp\left(\frac{a^{vv}_{11}}{\mu_1}x\right) \tilde v_s, \label{st1} \\\bar v_f =& \exp\left(\frac{a^{vv}_{22}}{\mu_2}x\right) \tilde v_f. \label{st4}
\end{align}
Then we obtain the scaled system:
\begin{align}
\notag\partial_t  \tilde w_s + v_s^\star \partial_x  \tilde w_s  =& \bar a^{ww}_{11}  \tilde w_s +\bar a^{ww}_{12} \tilde w_f \\&+ \bar a^{wv}_{11}(x) \bar v_s  +\bar a^{wv}_{12}(x) \bar v_{ f}, \label{bwv1}\\
\notag	\partial_t \tilde w_f + v_f^\star \partial_x \tilde w_f  =& \bar a^{ww}_{21} \tilde w_s +\bar a^{ww}_{22} \tilde w_f\\& + \bar a^{wv}_{21}(x) \bar v_{ s}   +\bar  a^{wv}_{22}(x) \bar v_{ f},\\
	\partial_t \bar v_{ s} - (\gamma p_s^\star - v_s^\star) \partial_x \bar v_{ s} =&\bar a^{vw}_{11}(x)\tilde w_s + \bar a^{vw}_{12}(x) \tilde w_{ f}   +\bar  a^{vv}_{12}(x) \bar v_{ f},\\
	\partial_t \bar v_{ f} - (\gamma p_f^\star-v_f^\star ) \partial_x  \bar v_{ f} =&\bar a^{vw}_{21}(x)\tilde w_s  +\bar a^{vw}_{22}(x)  \tilde w_f  +\bar a^{vv}_{21}(x)  \bar v_{ s},\\
	\tilde w_s(0,t) =& k_s \bar v_s(0,t),\\
    \tilde w_f(0,t) =& k_f \bar v_f(0,t),\\
	\bar v_{ s }(L,t) =&l_s U_s(t),\\
	\bar v_{ f }(L,t) =& l_f  U_f(t).\label{bwv2}
\end{align}
We denote the transports speeds as
\begin{align}
\epsilon_1 =& v_s^\star,\quad \epsilon_2 = v_f^\star,\\ \mu_1 =& ( \gamma p_s^\star-v_s^\star), \quad \mu_2 = ( \gamma p_f^\star-v_f^\star).
\end{align}
Note that the steady velocity of the fast lane is larger than that of the slow lane, the constant transport speeds satisfy the following inequalities,
\begin{align}
-\mu_1 < - \mu_2 < 0 < \epsilon_1<\epsilon_2.
\end{align}
where the constant coefficients $l_f$ and $l_s$ are defined as
\begin{align}
l_s = \exp\left(\frac{a^{vv}_{11}}{\mu_1}L\right), \quad  l_f = \exp\left(\frac{a^{vv}_{22}}{\mu_2}L\right).
\end{align}
The new in-domain coefficient matrix $\{\bar A\}$ is given by
\begin{align}
\bar A = \begin{bmatrix}
\bar  A^{ww} & \bar A^{wv}\\
\bar A^{vw} & \bar A^{vv}\\	
\end{bmatrix},
\end{align}
where the sub-matrices are obtained as
\begin{align}
\bar  A^{ww}=&A^{ww},\\
\bar A^{wv}(x) =& A^{wv} \begin{bmatrix}
\exp\left(-\frac{a^{vv}_{11}}{\mu_1}x\right)  &  0     \\
0  & \exp\left(-\frac{a^{vv}_{22}}{\mu_2}x\right) 
\end{bmatrix},\\
\bar A^{vw}(x) =& A^{vw} \begin{bmatrix}
\exp\left(\frac{a^{vv}_{11}}{\mu_1}x\right)  &  0     \\
0  & \exp\left(\frac{a^{vv}_{22}}{\mu_2}x\right)    
\end{bmatrix},\\
\bar  A^{vv}(x)=&A^{vv}\begin{bmatrix}
0 & \exp\left(\frac{a^{vv}_{22}}{\mu_2}x - \frac{a^{vv}_{11}}{\mu_1}x\right)\\
\exp\left(\frac{a^{vv}_{11}}{\mu_1}x- \frac{a^{vv}_{22}}{\mu_2}x\right)   
   & 0 
\end{bmatrix}.
\end{align}
Among the transformed sub-matrices, the elements of $\{\bar  A^{ww}\}$ are constant and the elements of $\{\bar A^{wv}(x)\}$, $\{\bar A^{vw}(x)\}$ and $\{\bar  A^{vv}(x)\}$ are spatially-varying coefficients. We summarize the transformation between $(\tilde w_s,\tilde w_f, \bar v_s, \bar v_f)$ and $(\tilde \rho_{ s },\tilde \rho_{ f }, \tilde v_s,\tilde v_f)$ from \eqref{rt2}, \eqref{rt1} and \eqref{st1}, \eqref{st4} as follows:
\begin{align}
\tilde\rho_{ s} = & \frac{\rho_s^\star}{\gamma p_s^\star}\left(\tilde w_s - \exp\left(-\frac{a^{vv}_{11}}{\mu_1}x\right) \bar v_s \right), \label{tc1}\\
\tilde \rho_{ f} = & \frac{\rho_f^\star}{\gamma p_f^\star}\left(\tilde w_f - \exp\left(-\frac{a^{vv}_{22}}{\mu_2}x\right) \bar v_s\right),\\
\tilde v_{ s}= &   \exp\left(-\frac{a^{vv}_{11}}{\mu_1}x\right) \bar v_s, \\
\tilde v_{ f}= & \exp\left(-\frac{a^{vv}_{22}}{\mu_2}x\right) \bar v_f. \label{tc4}
\end{align}
Then we introduce the backstepping transformation to the scaled $(\tilde w_i, \bar v_i)$-system in \eqref{bwv1}-\eqref{bwv2},
\begin{align}
\begin{bmatrix}
 \alpha_s(x,t)    \\
 \alpha_f(x,t)   
  \end{bmatrix} =& \begin{bmatrix}
  \tilde  w_s(x,t)      \\
  \tilde w_f(x,t)     
  \end{bmatrix},\label{pret1}\\
	\notag	\begin{bmatrix}
\beta_s(x,t)    \\
\beta_f(x,t)   
\end{bmatrix} =& \begin{bmatrix}
\bar v_s(x,t)      \\
\bar v_f(x,t)     
\end{bmatrix}-\int_{0}^{x} K(x,\xi)\begin{bmatrix}
\tilde  w_s(x,t)      \\
\tilde  w_f(x,t)     
\end{bmatrix}d\xi\\
&-\int_{0}^{x} L(x,\xi) \begin{bmatrix}
\bar v_s(x,t)      \\
\bar v_f(x,t)     
\end{bmatrix}d\xi, \label{pret2}
\end{align}
where the kernel matrices are denoted as
\begin{align}
K =& \begin{bmatrix}
K_{11}       & K_{12} \\
K_{21}       & K_{22}
\end{bmatrix},\quad L = \begin{bmatrix}
L_{11}       & L_{12} \\
L_{21}       & L_{22}
\end{bmatrix}.
\end{align}
The kernel variables $\{K\}$ and $\{L\}$ evolve in the triangular domain $\mathcal{T}=\{(x,\xi): 0 \leq \xi \leq x \leq 1\}$.
Taking derivative with respect to time and space on both sides of \eqref{pret1}-\eqref{pret2} along the solution of a target system given later, we obtain the following kernel equations. The kernels $\{K(x,\xi)\}$ and $\{L(x,\xi)\}$ are governed by
\begin{align}
&\mu_1 \partial_x K_{11}-\epsilon_1 \partial_{\xi} K_{11}=\bar a^{ww}_{11}K_{11} + \bar a^{ww}_{21}K_{12} + \bar a^{vw}_{11}L_{11}+\bar a^{vw}_{21} L_{12} ,\label{kn1}\\
&\mu_1 \partial_x K_{12}-\epsilon_2 \partial_{\xi} K_{12}=\bar a^{ww}_{12}K_{11} + \bar a^{ww}_{22}K_{12} + \bar a^{vw}_{12}L_{11}+\bar a^{vw}_{22} L_{12},\\
&\mu_2\partial_x K_{21}-\epsilon_1 \partial_{\xi} K_{21}=\bar a^{ww}_{11}K_{21} + \bar a^{ww}_{21}K_{22} + \bar a^{vw}_{11}L_{21}+\bar a^{vw}_{21}L_{22},\\
&\mu_2 \partial_x K_{22}-\epsilon_2 \partial_{\xi} K_{22}=\bar a^{ww}_{12}K_{21} + \bar a^{ww}_{22}K_{22} + \bar a^{vw}_{12}L_{21}+\bar a^{vw}_{22}L_{22},\\
&\mu_1\partial_x L_{11}+\mu_1 \partial_{\xi} L_{11}=  \bar a^{vv}_{21}L_{12} + \bar a^{vw}_{11}K_{11}+\bar a^{vw}_{21}K_{12},\\
&\mu_1\partial_x L_{12}+\mu_2 \partial_{\xi} L_{12}=\bar a^{vv}_{12}L_{11} + \bar a^{vw}_{12}K_{11}+\bar a^{vw}_{22}K_{12},\\
&\mu_2\partial_x L_{21}+\mu_1 \partial_{\xi} L_{21}=\bar a^{vv}_{21}L_{22} + \bar a^{vw}_{11}K_{21}+\bar a^{vw}_{21}K_{22},\\
&\mu_2\partial_x L_{22}+\mu_2 \partial_{\xi} L_{22}=\bar a^{vv}_{12}L_{21}+ \bar a^{vw}_{12}K_{21}+\bar a^{vw}_{22}K_{22},\\
&K_{11}(x,x)=-\frac{\bar a^{vw}_{11}(x)}{\epsilon_1 + \mu_1},\quad K_{12}(x,x)= -\frac{\bar a^{vw}_{12}(x)}{\epsilon_2+\mu_1},\\
&K_{21}(x,x)=-\frac{\bar a^{vw}_{21}(x)}{\epsilon_1 + \mu_2},\quad 
K_{22}(x,x)= - \frac{\bar a^{vw}_{22}(x)}{\epsilon_2 + \mu_2},\\
&L_{11}(x,0)=\frac{\epsilon_1 k_s}{\mu_1}K_{11}(x,0),\quad 
L_{12}(x,x)=-\frac{\bar a^{vv}_{12}(x)}{\mu_1-\mu_2} ,\\
&L_{12}(x,0)= \frac{\epsilon_2 k_f}{\mu_2}K_{12}(x,0),\quad 
L_{21}(x,x)= - \frac{\bar a^{vv}_{21}(x)}{\mu_2-\mu_1} ,\\
&L_{21}(L,\xi)= 0, \quad 
L_{22}(x,0)= \frac{\epsilon_2 k_f}{\mu_2}K_{22}(x,0).\label{kn16}
\end{align}
The well-possedness of the kernel equations \eqref{kn1}-\eqref{kn16} is proved using the method of characteristics and the successive approximations following the result for a general class of kernel system in~\cite{long}. There exists a unique solution $K,L \in L^\infty(\mathcal{T})$. Therefore, we establish the invertibility of the backstepping transformation \eqref{pret1},\eqref{pret2} and can study the stability of the following target system due to its equivalence to the $(\tilde w_i, \bar v_i)$-system. 

Note that we impose an artificial boundary condition $L_{21}(L,\xi)$ in \eqref{kn16} for the well-posedness of the kernel equations. This leads to one degree of freedom in backstepping transformation of the hyperbolic system as well as the following control design. The stabilization of the following target system is achieved with two controllers and the one degree of freedom enables the coordination between the two VSLs.

With the backstepping transformation and the above kernel equations, we map the  $(\bar w_i, \bar v_i)$-system to the cascade target system $(\alpha_i, \beta_i)$,
\begin{align}
	\notag	\partial_t \alpha_s + v_s^\star \partial_x \alpha_s  =& \bar a^{ww}_{11} \alpha_s +\bar a^{ww}_{12}\alpha_f \\\notag&+ \bar a^{wv}_{11}(x) \beta_s  +\bar a^{wv}_{12}(x) \beta_{ f} \\\notag&+\int_{0}^{x}b_{11}(x,\xi)\alpha_{ s}(\xi)d\xi\\\notag&+\int_{0}^{x}b_{12}(x,\xi)\alpha_{f}(\xi)d\xi\\\notag	&+\int_{0}^{x}c_{11}(x,\xi)\beta_{ s}(\xi)d\xi\\&+\int_{0}^{x}c_{12}(x,\xi)\beta_{ f}(\xi)d\xi,\label{t1}\\
	\notag	\partial_t \alpha_f + v_f^\star \partial_x \alpha_f  =& \bar a^{ww}_{21} \alpha_s +\bar a^{ww}_{22}\alpha_f \\\notag&+ \bar a^{wv}_{21}(x) \beta_{ s}   +\bar  a^{wv}_{22}(x) \beta_{ f}\\\notag &+\int_{0}^{x}b_{21}(x,\xi)\alpha_{ s}(\xi)d\xi\\\notag&+\int_{0}^{x}b_{22}(x,\xi)\alpha_{f}(\xi)d\xi \\	
\notag	&+\int_{0}^{x}c_{21}(x,\xi)\beta_{ s}(\xi)d\xi\\&+\int_{0}^{x}c_{22}(x,\xi)\beta_{ f}(\xi)d\xi,\\
	\partial_t \beta_{ s} - (\gamma p_s^\star - v_s^\star) \partial_x \beta_{ s} =&0, \label{c1}\\
	\partial_t \beta_{ f} - (\gamma p_f^\star-v_f^\star ) \partial_x  \beta_{ f} =&\theta(x)  \beta_{ s}(0,t)  ,\label{c2} \\
    \alpha_f(0,t) =& k_f \beta_f(0,t),\\
	\alpha_s(0,t) =& k_s \beta_s(0,t),\\
	\beta_{ s }(L,t) =&0,\label{t2}\\
	\beta_{ f }(L,t) =& 0,\label{beta0}
\end{align}
where the spatially varying parameter matrices $\{B\}$ and $\{C\}$ are denoted as
\begin{align}
B = \begin{bmatrix}
b_{11} & b_{12}\\
b_{21} & b_{22}\\	
\end{bmatrix},\quad
C = \begin{bmatrix}
c_{11} & c_{12}\\
c_{21} & c_{22}\\	
\end{bmatrix},
\end{align} 
 and given by the following equations in the matrix form,
\begin{align}
B(x,\xi)  =&\bar A^{wv} K(x,\xi)+\int_{\xi}^{x} B(x,s)K(s,\xi)d\xi,\\
C(x,\xi) = & \bar A^{wv}L(x,\xi)+\int_{\xi}^{x}C(x,s)L(s,\xi)d\xi,
\end{align}
and $\theta(x)$ is obtained from the kernel variables $K_{21}$ and $ L_{21}$,
\begin{align}
	\theta(x)= -\epsilon_1 k_s K_{21}(x,0)-\mu_1 L_{21}(x,0).
\end{align}
Considering the cascade structure of the target system, the following conclusion is arrived.
\begin{lem}
Consider the target system \eqref{t1}-\eqref{c2} and actuated boundary conditions \eqref{t1}-\eqref{beta0}, the zero equilibrium
\begin{align}
\alpha_f(x,t) = \alpha_s (x,t) = \beta_{ f}(x,t) =\beta_{ s}(s,t) \equiv 0,
\end{align}
is reached in finite time $t=t_f$, where
\begin{align}
t_f = \frac{L}{v_s^\star} +\frac{L}{\gamma p_f^\star - v_f^\star} + \frac{L}{\gamma p_s^\star - v_s^\star}. \label{tf}
\end{align}
\end{lem}
\begin{pf*}{Proof.}
By solving \eqref{c2} and \eqref{beta0} directly, we obtain that after $t>\frac{L}{\gamma p_f^\star - v_f^\star}$, 
\begin{align}
\beta_{ f}(x,t) \equiv 0.
\end{align}
Using the cascade structure of the target system in \eqref{c1}, we have after $t>\frac{L}{\gamma p_f^\star - v_f^\star}+\frac{L}{\gamma p_s^\star - v_s^\star}$,
\begin{align}
\beta_{ s}(s,t) \equiv 0.
\end{align}
Then after $t>\frac{L}{\gamma p_f^\star - v_f^\star}+\frac{L}{\gamma p_s^\star - v_s^\star}+\frac{L}{v_s^\star}$, we obtain that 
\begin{align}
\alpha_f(x,t) \equiv 0,\quad \alpha_s (x,t)\equiv 0,
\end{align}
which concludes the proof. 
\end{pf*}
Boundary conditions \eqref{t2},\eqref{beta0} and backstepping transformation \eqref{pret2} yield full-state feedback control laws given by $(\bar w_i,\bar v_i)$,
\begin{align}
\notag	\begin{bmatrix}
	l_s U_s(t)      \\
	l_f U_f(t)     
	\end{bmatrix}=&\int_{0}^{L} \begin{bmatrix}
	K_{11}(L,\xi)       & K_{12}(L,\xi) \\
	K_{21} (L,\xi)      & K_{22}(L,\xi)
	\end{bmatrix}\begin{bmatrix}
	\bar w_s(x,t)      \\
	\bar w_f(x,t)     
	\end{bmatrix}d\xi\\
	&+\int_{0}^{L} \begin{bmatrix}
	L_{11} (L,\xi)      & L_{12}(L,\xi) \\
	L_{21}  (L,\xi)     & L_{22}(L,\xi)
	\end{bmatrix}\begin{bmatrix}
	\bar v_s(x,t)      \\
	\bar v_f(x,t)     
	\end{bmatrix}d\xi.  \label{full_control}
\end{align}
Using the invertible transformation \eqref{tc1}-\eqref{tc4}, we obtain the full-state feedback control laws given in traffic flow variables $(\rho_{ f}, v_{ f}, \rho_{ s}, v_{ s})$ and the steady states $(\rho_{ f}^\star, v_{ f}^\star, \rho_{ s}^\star, v_{ s}^\star)$.
We reach the main stabilization result of full-state feedback control design.
\begin{thm}
	Consider the two-lane traffic ARZ model in \eqref{rv1}-\eqref{rv4} with boundary conditions \eqref{bc1}-\eqref{bc4}, initial conditions $\rho_{ f}(x,0), v_{ f}(x,0), \rho_{ s}(x,0), v_{ s}(x,0) \in {L}^\infty\left([0,L]\right)$ and the following control laws	
\begin{align}
		\notag U_s(t)=&\exp\left(-\frac{a^{vv}_{11}}{\mu_1}L\right)  \int_{0}^{L} \frac{\gamma p_s^\star}{\rho_s^\star} K_{11}(L,\xi) \left(\rho_{ s}(\xi, t)-\rho_{ s}^\star\right) 
		\\ \notag  &\!\!\!\!\!+ \frac{\gamma p_f^\star}{\rho_f^\star}K_{12}(L,\xi) \left(\rho_{ f}(\xi, t)-\rho_{ f}^\star\right)
		\\ \notag &\!\!\!\!\! +\!\left[\!\!K_{11}(L,\xi)+L_{11}(L,\xi)\exp\left(\!\!\frac{a_{11}^{vv}}{\mu_1}\xi\!\right)\!\!\right]\!\! \left(v_{ s}(\xi, t)-v_{ s}^\star\right) \\
		&\!\!\!\!\!+\!\left[\!\!K_{12}(L,\xi)+L_{12}(L,\xi)\exp\left(\!\!\frac{a_{22}^{vv}}{\mu_2}\xi\!\right)\!\!\right]\!\! \left(v_{ f}(\xi, t)-v_{ f}^\star\right)\!d\xi, \label{Us}\\
		\notag U_f(t)=&\exp\left(-\frac{a^{vv}_{22}}{\mu_2}L\right)\int_{0}^{L} \frac{\gamma p_s^\star}{\rho_s^\star}  K_{21}(L,\xi)  \left(\rho_{ s}(\xi, t)-\rho_{ s}^\star\right) 
		\\ \notag  &\!\!\!\!\!+\frac{\gamma p_f^\star}{\rho_f^\star}  K_{22}(L,\xi) \left(\rho_{ f}(\xi, t)-\rho_{ f}^\star\right) 
		\\ \notag  &\!\!\!\!\!+\!\left[\!\!K_{21}(L,\xi)+L_{21}(L,\xi)\exp\left(\!\!\frac{a_{11}^{vv}}{\mu_1}\xi\!\right)\!\!\right]\!\! \left(v_{ s}(\xi, t)-v_{ s}^\star\right)\\
		&\!\!\!\!\!+\!\left[\!\!K_{22}(L,\xi) + L_{22}(L,\xi)\exp\left(\!\!\frac{a_{22}^{vv}}{\mu_2}\xi\!\right)\!\!\right]\!\! \left(v_{ f}(\xi, t)-v_{ f}^\star\right)\!d\xi,\label{Uf}
\end{align}
where the kernels are obtained by solving \eqref{kn1}-\eqref{kn16}. The steady states $(\rho_{ f}^\star, v_{ f}^\star, \rho_{ s}^\star, v_{ s}^\star)$ are finite-time stable and the convergence is reached in $t_f$ given in \eqref{tf}.
\end{thm}
\begin{pf*}{Proof.}
	Lemma 1 for the closed-loop target system in \eqref{t1}-\eqref{beta0} with the existence of the backstepping transformation in \eqref{pret1},\eqref{pret2} yields the convergence of the states variables $(\tilde w_s$,$\tilde w_f$,$\bar v_s$,$\bar v_f)$ defined by \eqref{bwv1}-\eqref{bwv2} to zero for $t>t_f$. Given the transformation in \eqref{tc1}-\eqref{tc4}, the finite-time convergence to zero is arrived for the linearized state variables $(\tilde \rho_s(x,t), \tilde \rho_f(x,t), \tilde v_s(x,t), \tilde v_f(x,t))$, which yields the convergence of the two-lane ARZ PDE model by $(\rho_s(x,t), \rho_f(x,t), v_s(x,t), v_f(x,t))$ to the steady states.
\end{pf*}
\section{Collocated Observer Design and Output Feedback Control}
In this section, we develop a collocated observer by taking measurement of density states at the outlet of the segment,
\begin{align}
	y_s(t) =& \tilde \rho_s(L,t),\\
	y_f(t) =& \tilde \rho_f(L,t). 
\end{align}
Using the state estimates obtained form the observer design and the full-state feedback control laws, we construct output feedback controllers.

Note that the anti-collocated observer can also be designed here by taking measurement of velocity states $\tilde v_{ s}(0,t)$ and  $\tilde v_{f}(0,t)$ at the inlet. The anti-collocated observer design is trivial in our case which presents as a copy of the $(\tilde w_f, \tilde w_s, \bar v_s, \bar v_f)$-system. More importantly, collocated observer design is practical in implementation along with the full-state feedback control design. 

\subsection{Collocated observer design}
For state estimation of the scaled system in \eqref{bwv1}-\eqref{bwv2}, we obtain the measurement of $\tilde w_s(L,t)$ and $\tilde w_f(L,t)$ from \eqref{rt2},\eqref{rt1},
\begin{align}
	Y_s(t)=& \tilde w_s(L,t)  = \frac{\gamma p_s^\star}{\rho_s^\star} \tilde \rho_{ s}(L,t)  + \tilde v_{ s}(L,t),\\
	Y_f(t)=& \tilde w_f(L,t)  = \frac{\gamma p_f^\star}{\rho_f^\star} \tilde \rho_{ f}(L,t)  + \tilde v_{ f}(L,t),
\end{align}
thus the values of $Y_s(t)$ and $Y_f(t)$ are obtained from $y_s(t)$, $y_f(t)$ and control inputs $U_s(t)$, $U_f(t)$,
\begin{align}
Y_s(t)=& \frac{\gamma p_s^\star}{\rho_s^\star} y_s(t)  +  U_s(t),\\
Y_f(t)=& \frac{\gamma p_f^\star}{\rho_f^\star} y_f(t)  +  U_f(t).	
\end{align}
The observer equations $(\hat w_f, \hat w_s, \hat u_s, \hat u_f)$ that estimate $(\tilde w_f, \tilde w_s, \bar v_s, \bar v_f)$ read as follows:
\begin{align}
\notag	\partial_t \hat w_s + v_s^\star \partial_x \hat w_s  =&  \bar a^{ww}_{11}  \hat w_s +\bar a^{ww}_{12} \hat w_f \\\notag &+ \bar a^{wv}_{11}(x) \hat u_s  +\bar a^{wv}_{12}(x) \hat u_{ f}\\
& + p_{11}(x)\check w_s(L,t) + p_{12}(x)\check w_f(L,t), \label{ob1}\\
\notag	\partial_t \hat w_f + v_f^\star \partial_x \hat w_f  =& \bar a^{ww}_{21} \hat w_s +\bar a^{ww}_{22} \hat w_f \\\notag &+ \bar a^{wv}_{21}(x) \hat u_{ s}   +\bar  a^{wv}_{22}(x) \hat u_{ f}\\
	& + p_{21}(x)\check w_s(L,t) + p_{22}(x)\check w_f(L,t), \\
\notag	\partial_t \hat u_{ s} - (\gamma p_s^\star - v_s^\star) \partial_x \hat u_{ s} =&\bar a^{vw}_{11}(x)\hat w_s + \bar a^{vw}_{12}(x) \hat w_{ f}   +\bar  a^{vv}_{12}(x) \hat u_{ f}\\
	& + q_{11}(x)\check w_s(L,t) + q_{12}(x)\check w_f(L,t), \\
\notag	\partial_t \hat u_{ f} - (\gamma p_f^\star-v_f^\star ) \partial_x  \hat u_{ f} =&\bar a^{vw}_{21}(x)\hat w_s  +\bar a^{vw}_{22}(x)  \hat w_f  +\bar a^{vv}_{21}(x)  \hat u_{ s}\\
	& + q_{21}(x)\check w_s(L,t) + q_{22}(x)\check w_f(L,t), \label{ob4} \\
	\hat w_s(0,t) =& k_s \hat u_{ s }(0,t),\\
	\hat w_f(0,t) =& k_f \hat u_{ f }(0,t),\\
	\hat u_{ s }(L,t) =&l_s U_s(t),\\
	\hat u_{ f }(L,t) =& l_f U_f(t).\label{obl}
\end{align}
The output injections in \eqref{ob1}-\eqref{ob4} are defined as
\begin{align}
\check w_s(L,t) =& \tilde w_s(L,t) -\hat w_s(L,t),\\
\check w_f(L,t) =& \tilde w_f(L,t) -\hat w_f(L,t).
\end{align}
The observer output injection gains matrices $\{P\}$ and $\{Q\}$ are denoted as
\begin{align}
P = \begin{bmatrix}
p_{11}       & p_{12} \\
p_{21}       & p_{22}
\end{bmatrix},\quad	Q = \begin{bmatrix}
q_{11}       & q_{12} \\
q_{21}       & q_{22}
\end{bmatrix}.
\end{align}
The output injection gains are to be designed so that the output injection terms can drive the estimation error system of the observer to converge to zero in finite-time. 

The estimation errors are defined as 
\begin{align}
\check w_s =& \tilde w_s - \hat w_s, \quad \check v_{ s} =\bar v_{ s} - \hat u_{ s},\\
\check w_f =&\tilde w_f - \hat w_f, \quad \check v_{ f} =\bar v_{ f} - \hat u_{ f}.		
\end{align}
The error system $(\check w_s, \check w_f, \check v_f, \check v_s)$ of the observer is given by subtracting \eqref{ob1}-\eqref{obl} from \eqref{bwv1}-\eqref{bwv2},
\begin{align}
	\notag	\partial_t  \check w_s + v_s^\star \partial_x  \check w_s  =& \bar a^{ww}_{11}  \check w_s +\bar a^{ww}_{12} \check w_f \\\notag&+ \bar a^{wv}_{11}(x) \check v_s  +\bar a^{wv}_{12}(x) \check v_{ f}\\
	& - p_{11}(x)\check w_s(L,t) - p_{12}(x)\check w_f(L,t), \label{err1}\\
	\notag	\partial_t \check w_f + v_f^\star \partial_x \check w_f  =& \bar a^{ww}_{21} \check w_s +\bar a^{ww}_{22} \check w_f \\\notag& + \bar a^{wv}_{21}(x) \check v_{ s}   +\bar  a^{wv}_{22}(x) \check v_{ f},\\
	& - p_{21}(x) \check w_s(L,t) - p_{22}(x)\check w_f(L,t), \\
	\notag	\partial_t \check v_{ s} - (\gamma p_s^\star - v_s^\star) \partial_x \check v_{ s} =&\bar a^{vw}_{11}(x)\check w_s + \bar a^{vw}_{12}(x) \check w_{ f}  + \bar  a^{vv}_{12}(x) \check v_{ f}\\
	& - q_{11}(x) \check w_s(L,t) - q_{12}(x)\check w_f(L,t), \\
	\notag	\partial_t \check v_{ f} - (\gamma p_f^\star-v_f^\star ) \partial_x  \check v_{ f} =&\bar a^{vw}_{21}(x)\check w_s  +\bar a^{vw}_{22}(x)  \check w_f  +\bar a^{vv}_{21}(x)  \check v_{ s},\\
	& - q_{21}(x)\check w_s(L,t) - q_{22}(x)\check w_f(L,t), \\
	\check w_s(0,t) =& k_s \check v_s(0,t),\\
	\check w_f(0,t) =& k_f \check v_f(0,t),\\
	\check v_{ s }(L,t) =&0,\\
	\check v_{ f }(L,t) =&0,\label{errl}
	\end{align}
We apply the backstepping transformation to the error system given by
\begin{align}
\begin{bmatrix}
\check w_s(x,t)    \\
\check w_f(x,t)   
\end{bmatrix} =& \begin{bmatrix}
\check \alpha_s(x,t)      \\
\check \alpha_f(x,t)     
\end{bmatrix}+\int_{x}^{L} M(x,\xi) \begin{bmatrix}
\check \alpha_s(\xi,t)      \\
\check \alpha_f(\xi,t)     
\end{bmatrix} d\xi, \label{bkst1} \\
\begin{bmatrix}
\check v_s(x,t)      \\
\check v_f(x,t)     
\end{bmatrix}	 =& 
\begin{bmatrix}
\check\beta_s(x,t)    \\
\check \beta_f(x,t)   
\end{bmatrix}
+\int_{x}^{L}N(x,\xi)
\begin{bmatrix}
\check \alpha_s(\xi,t)      \\
\check \alpha_f(\xi,t)     
\end{bmatrix}d\xi, \label{bkst2}
\end{align}
where the kernel matrices $\{M\}$, $\{N\}$ are denoted as
\begin{align}
M = \begin{bmatrix}
M_{11}       & M_{12} \\
M_{21}       & M_{22}
\end{bmatrix},\quad
N = \begin{bmatrix}
N_{11}       & N_{12} \\
N_{21}       & N_{22}
\end{bmatrix}.
\end{align}
The kernels $\{M\}$, $\{N\}$ evolve in the triangular domain $\mathcal{T}=\{(x,\xi): 0  \leq x \leq \xi \leq L\}$ and are defined later.
We map the error system in \eqref{err1}-\eqref{errl} into the following cascade target system
\begin{align}
\notag	\partial_t \check \alpha_s + v_s^\star \partial_x \check \alpha_s  =&\bar a_{11}^{ww} \check\alpha_s + \bar a_{11}^{wv}(x) \check \beta_s  + \bar a_{12}^{wv}(x) \check \beta_f \\\notag &+\int_{x}^{L}d_{11}(x,\xi)\check\beta_{ s}(\xi)d\xi\\
&+\int_{x}^{L}d_{12}(x,\xi)\check\beta_{f}(\xi)d\xi, \label{1ot}\\
\notag	\partial_t \check \alpha_f + v_f^\star \partial_x  \check \alpha_f  =&\bar a_{22}^{ww} \check\alpha_f +  \bar a_{21}^{wv} (x) \check \beta_{ s} +\bar a_{22}^{wv}(x) \check\beta_{ f}  \\\notag&+\int_{x}^{L}d_{21}(x,\xi)\check \beta_{ s}(\xi)d\xi\\&+\int_{x}^{L}d_{22}(x,\xi) \check \beta_{f}(\xi)d\xi,\label{2ot}\\
\notag\partial_t \check \beta_{ s} - (\gamma p_s^\star - v_s^\star) \partial_x \check \beta_{ s} =&\bar a_{12}^{vv}(x) \check \beta_f \\\notag &+\int_{x}^{L}f_{11}(x,\xi)\check\beta_{ s}(\xi)d\xi\\&+\int_{x}^{L}f_{12}(x,\xi)\check\beta_{f}(\xi)d\xi,\label{3ot}\\
\notag\partial_t \check \beta_{ f} - (\gamma p_f^\star-v_f^\star ) \partial_x  \check \beta_{ f} =&\bar a_{21}^{vv}(x) \check \beta_s \\\notag &+\int_{x}^{L}f_{21}(x,\xi)\check\beta_{ s}(\xi)d\xi\\&+\int_{x}^{L}f_{22}(x,\xi)\check\beta_{f}(\xi)d\xi,\label{4ot}\\
 \check \alpha_s(0,t) =& k_s \check \beta_s(0,t),\label{1ob}\\
\check \alpha_f(0,t) =& k_f \check \beta_f(0,t) - \int_{0}^{L}\lambda(x) \check\alpha_s(x,t) d\xi ,\\
\check \beta_{ s }(L,t) =&0,\label{3ob} \\
\check \beta_{ f }(L,t) =& 0,\label{4ob}
\end{align}
where the coefficient matrices $\{D\}$ and $\{F\}$ are 
denoted as
\begin{align}
D = \begin{bmatrix}
d_{11}       & d_{12} \\
d_{21}       & d_{22}
\end{bmatrix},\quad
F = \begin{bmatrix}
f_{11}       & f_{12} \\
f_{21}       & f_{22}
\end{bmatrix}.
\end{align}
and given by
\begin{align}
D(x,\xi)  =&-M(x,\xi)\bar A^{wv} +\int_{\xi}^{x} M(x,s)D(s,\xi)d\xi,\\
F(x,\xi) = &-N(x,\xi)\bar A^{ww}+\int_{\xi}^{x}N(x,s)F(s,\xi)d\xi.
\end{align}
The spatially varying coefficient $\lambda(x)$ is obtained from the kernel variables
\begin{align}
	\lambda(x) = M_{21}(0,x) - k_f N_{21} (0,x).
\end{align}

\begin{lem}
	Consider the target system \eqref{1ot}-\eqref{4ot} with the boundary conditions \eqref{1ob}-\eqref{4ob}. The zero equilibrium is reached in finite time $t=t_o$,
	where 
	\begin{align}
	t_o = \frac{L}{v_s^\star} + \frac{L}{v_f^\star} + \frac{L}{\gamma p_f^\star - v_f^\star}. \label{to}
	\end{align}
\end{lem}
\vspace*{-1cm}
\begin{pf*}{Proof.}
	Noting the cascade structure of $\check\alpha_s$, $\check\alpha_f$  and $\check\beta_s$, $\check\beta_f$ system, $\check\beta$ variables appear as the right hand source terms in $\check\alpha$ equations and through the inlet boundaries. The integral of variable $\check\alpha_s$ enters the boundary condition of $\check\alpha_f$. Therefore we solve the target system explicitly by recursion. The $\check\beta$-system is independent of the $\check\alpha$ system. Given the boundary conditions in \eqref{3ob}, \eqref{4ob}, the explicit solutions hold for $\check\beta_s$ and $\check\beta_f$ after $t> \frac{1}{\gamma p_f^\star - v_f^\star}$, 
	\begin{align}
		\check\beta_{ s}(x,t) \equiv 0, \quad \check\beta_{f}(x,t) \equiv 0.
	\end{align}
	When $t> \frac{1}{\gamma p_f^\star - v_f^\star}$, $\check\alpha$-system becomes
	\begin{align}
		\partial_t \check \alpha_s + v_s^\star \partial_x \check \alpha_s  =&\bar a_{11}^{ww} \check\alpha_s,\\
		\partial_t \check \alpha_f + v_f^\star \partial_x  \check \alpha_f  =&\bar a_{22}^{ww} \check\alpha_f,\\
		 \check \alpha_s(0,t) =& 0,\\
		\check \alpha_f(0,t) =& - \int_{0}^{L}\lambda(x) \check\alpha_s(x,t) d\xi.
	\end{align}
    After $t > \frac{1}{v_s^\star}+\frac{1}{\gamma p_f^\star - v_f^\star}$, we have $\check\alpha_s$ satisfies  
    \begin{align}
    \check\alpha_{ s}(x,t) \equiv 0.
    \end{align}
    Then $ \check\alpha_{f}(x,t) \equiv 0$ follows after another time period $\frac{1}{v_f^\star}$. Therefore, $\check\alpha$-system eventually identically vanishes for 
    \begin{align}
    t_o = \frac{1}{v_s^\star} + \frac{1}{v_f^\star} + \frac{1}{\gamma p_f^\star - v_f^\star},
    \end{align}
    which concludes the proof.
\end{pf*}
Taking spatial and temporal derivatives of the backstepping transformation \eqref{bkst1},\eqref{bkst2} along the target system \eqref{1ot}-\eqref{4ob}, then plugging into the error system \eqref{err1}-\eqref{errl}, we obtain the kernel equations that govern 
the kernels $\{M(x,\xi)\}$ and $\{N(x,\xi)\}$,
\begin{align}
	\epsilon_1 \partial_x M_{11} + \epsilon_1 \partial_{\xi} M_{11} =& - \bar a^{ww}_{12}M_{21} - \bar a^{wv}_{11}N_{11}-\bar a^{wv}_{12} N_{21} ,\label{mn1}\\
\notag	\epsilon_1 \partial_x M_{12} + \epsilon_2 \partial_{\xi} M_{12} =& -\bar a^{ww}_{11}M_{12} - \bar a^{ww}_{12}M_{22} \\ &- \bar a^{wv}_{11}N_{12}-\bar a^{wv}_{12} N_{22},\\
\notag	\epsilon_2\partial_x M_{21} + \epsilon_1\partial_{\xi} M_{21} =& -\bar a^{ww}_{21}M_{11} - \bar a^{ww}_{22}M_{21} \\ &- \bar a^{wv}_{21}N_{11} - \bar a^{wv}_{22}N_{21},\\
	\epsilon_2 \partial_x M_{22} + \epsilon_2 \partial_{\xi} M_{22}=& -\bar a^{ww}_{21}M_{12}  - \bar a^{wv}_{21}N_{12}-\bar a^{wv}_{22} N_{22},\\
	\mu_1\partial_x N_{11}-\epsilon_1 \partial_{\xi} N_{11} =& \bar a^{ww}_{11} N_{11} +  \bar a^{vv}_{12}N_{21} + \bar a^{vw}_{11}M_{11}+\bar a^{vw}_{12}M_{21},\\
	\mu_1\partial_x N_{12}-\epsilon_2 \partial_{\xi} N_{12} =& \bar a^{ww}_{22} N_{12} + \bar a^{vv}_{12}N_{22} + \bar a^{vw}_{11}M_{12}+\bar a^{vw}_{12}M_{22},\\
	\mu_2\partial_x N_{21}-\epsilon_1\partial_{\xi} N_{21} =& \bar a^{ww}_{11} N_{21} + \bar a^{vv}_{21}N_{11} + \bar a^{vw}_{21}M_{11}+\bar a^{vw}_{22}M_{21},\\
	\mu_2\partial_x N_{22}-\epsilon_2\partial_{\xi} N_{22} =& \bar a^{ww}_{22} N_{22} + \bar a^{vv}_{21}N_{12}+ \bar a^{vw}_{21}M_{12}+\bar a^{vw}_{22}M_{22},\label{mn16}
\end{align}
\vspace*{-0.8cm}
\begin{align}
	&N_{11}(x,x)=\frac{\bar a^{vw}_{11}(x)}{\epsilon_1 + \mu_1},
	\quad N_{12}(x,x)= \frac{\bar a^{vw}_{12}(x)}{\epsilon_2+\mu_1},\\
	&N_{21}(x,x)=\frac{\bar a^{vw}_{21}(x)}{\epsilon_1 + \mu_2},\quad N_{22}(x,x)=  \frac{\bar a^{vw}_{22}(x)}{\epsilon_2 + \mu_2},\\
	&M_{11}(0,\xi)= k_s N_{11}(0,\xi),\quad M_{22}(0,\xi)=  k_f N_{22}(0,\xi),\\
	&M_{12}(0,\xi)=  k_s N_{12}(0,\xi),\quad M_{21}(x,x)= - \frac{\bar a^{ww}_{21}}{\epsilon_2-\epsilon_1} ,\\
	&M_{21}(x,L)= - \frac{\bar a^{ww}_{21}}{\epsilon_2-\epsilon_1} , \quad M_{12}(x,x)=-\frac{\bar a^{ww}_{12}}{\epsilon_1-\epsilon_2}.
	\label{mn2}
\end{align}

Considering the following variables by defining 
\begin{align}
	\bar x = L-x,\quad \bar \xi = L-\xi,
\end{align}
and new kernels $\bar M\left(\bar x,\bar \xi\right)$ and $\bar N\left(\bar x,\bar \xi\right) $ 
\begin{align}
	\bar M\left(\bar x,\bar \xi\right) =& M(L-\bar x,L-\bar \xi)  = M(x,\xi),\\
	\bar N\left(\bar x,\bar \xi\right) =& N(L-\bar x,L-\bar \xi)  = N(x,\xi),
\end{align}
which are defined in the triangular domain $\mathcal{D}=\{(\bar x,\bar \xi): 0  \leq \bar\xi \leq \bar x \leq L\}$. We find that the following kernel equations obtained from \eqref{mn1} and \eqref{mn2} have the same structure with the controller kernel system in \eqref{kn1}-\eqref{kn16}. 
	\begin{align}
	\epsilon_1 \partial_{\bar x} \bar M_{11} + \epsilon_1 \partial_{\bar \xi} \bar M_{11} =&  \bar a^{ww}_{12} \bar M_{21} + \bar a^{wv}_{11}\bar N_{11} +\bar a^{wv}_{12} \bar N_{21} , \\
\notag	\epsilon_1 \partial_{\bar x}  \bar M_{12} + \epsilon_2 \partial_{\bar \xi} \bar M_{12}=& \bar a^{ww}_{11}\bar M_{12} + \bar a^{ww}_{12}\bar M_{22} \\&+ \bar a^{wv}_{11}\bar N_{12}+\bar a^{wv}_{12}\bar N_{22},\\
\notag	\epsilon_2\partial_{\bar x}  \bar M_{21} + \epsilon_1\partial_{\bar \xi} \bar M_{21}=& \bar a^{ww}_{21}\bar M_{11} + \bar a^{ww}_{22}\bar M_{21} \\&+  \bar a^{wv}_{21}\bar N_{11} + \bar a^{wv}_{22}\bar N_{21},\\
	\epsilon_2 \partial_{\bar x}  \bar M_{22} + \epsilon_2 \partial_{\bar \xi} \bar M_{22}=&\bar a^{ww}_{21}\bar M_{12} + \bar a^{wv}_{21}\bar N_{12} + \bar a^{wv}_{22}\bar N_{22},\\
\notag	\mu_1 \partial_{\bar x}  \bar N_{11} - \epsilon_1 \partial_{\bar \xi}\bar N_{11}=& -\bar a^{ww}_{11}\bar N_{11} -  \bar a^{vv}_{12}\bar N_{21} \\&-\bar a^{vw}_{11}\bar M_{11} - \bar a^{vw}_{12}\bar M_{21},\\
\notag	\mu_1 \partial_{\bar x}  \bar N_{12} - \epsilon_2 \partial_{\bar \xi}\bar N_{12}=&-\bar a^{ww}_{22} \bar N_{12} - \bar a^{vv}_{12}\bar N_{22} \\&- \bar a^{vw}_{11}\bar M_{12} - \bar a^{vw}_{12}\bar M_{22},\\
\notag	 \mu_2\partial_{\bar x}  \bar N_{21} - \epsilon_1 \partial_{\bar \xi} \bar N_{21}=&-\bar a^{ww}_{11} \bar N_{21} - \bar a^{vv}_{21} \bar N_{11} \\&- \bar a^{vw}_{21}\bar M_{11} - \bar a^{vw}_{22}\bar M_{21},\\
\notag	\mu_2 \partial_{\bar x}  \bar N_{22} - \epsilon_2 \partial_{\bar \xi} N_{22}=&-\bar a^{ww}_{22} \bar N_{22} - \bar a^{vv}_{21} \bar N_{12} \\&- \bar a^{vw}_{21}\bar M_{12}-\bar a^{vw}_{22}\bar M_{22},
\end{align}
\vspace*{-0.8cm}
\begin{align}
	\bar N_{11}(\bar x,\bar x)=&\frac{\bar a^{vw}_{11}(L-\bar x)}{\epsilon_1 + \mu_1},
	\quad \bar N_{12}(\bar x,\bar x)= \frac{\bar a^{vw}_{12}(L-\bar x)}{\epsilon_2+\mu_1},\\
	\bar N_{21}(\bar x,\bar x)=&\frac{\bar a^{vw}_{21}(L-\bar x)}{\epsilon_1 + \mu_2},\quad 
	\bar N_{22}(\bar x,\bar x)=  \frac{\bar a^{vw}_{22}(L-\bar x)}{\epsilon_2 + \mu_2},\\
	\bar M_{11}(L,\bar \xi)= &k_s \bar N_{11}(L,\bar \xi),\quad\bar M_{22}(L,\bar \xi)=  k_f \bar N_{22}(L,\bar \xi),\\
	\bar M_{12}(L,\bar \xi)=  &k_s \bar N_{12}(L,\bar \xi),\quad
	\bar M_{21}(\bar x,\bar x)= - \frac{\bar a^{ww}_{21}}{\epsilon_2-\epsilon_1} ,\\
	\bar M_{21}(\bar x,0)=& 0, \quad \bar M_{12}(\bar x,\bar x)=-\frac{\bar a^{ww}_{12}}{\epsilon_1-\epsilon_2}.
	\end{align}

The well-posedness of the above kernel system is obtained following the same steps of the proof for \eqref{kn1}-\eqref{kn16}. Therefore, there exists a unique solution $M,N\in L^\infty(\mathcal{T})$. The stability of target system \eqref{1ot}-\eqref{4ob} is equivalent to the error system \eqref{err1}-\eqref{errl}. The artificial boundary condition $\bar M_{21}(\bar x,0)$ in \eqref{kn16} is imposed for the well-posedness of the kernel equations. 

The observer gains matrices $\{P(x)\}$ and $\{Q(x)\}$ are obtained from the kernel matrices
\begin{align}
		P(x) =   M(x,L)\begin{bmatrix}
		v_s^\star      & 0 \\
		0       & v_f^\star
		\end{bmatrix},\label{PM}\\
	    Q(x) =   N(x,L)\begin{bmatrix}
	    v_s^\star      & 0 \\
	    0       & v_f^\star
	    \end{bmatrix}.\label{QN}
\end{align}
Note that the states estimation of the original traffic flow variables $(\hat\rho_{ f}, \hat v_{ f}, \hat \rho_{ s}, \hat v_{ s})$ are obtained by the invertible transformation given in the following,
\begin{align}
\hat \rho_{ s} = & \rho_{ s}^\star + \frac{\rho_s^\star}{\gamma p_s^\star}\left(\hat w_s - \exp\left(-\frac{a^{vv}_{11}}{\mu_1}x\right) \hat u_s \right), \label{ot1}\\
\hat \rho_{ f} = & \rho_{ f}^\star + \frac{\rho_f^\star}{\gamma p_f^\star}\left(\hat w_f - \exp\left(-\frac{a^{vv}_{22}}{\mu_2}x\right) \hat u_f \right),\\
\hat v_{ s}= &  v_{ s}^\star + \exp\left(-\frac{a^{vv}_{11}}{\mu_1}x\right) \hat u_s, \\
\hat v_{ f}= & v_{ f}^\star + \exp\left(-\frac{a^{vv}_{22}}{\mu_2}x\right) \hat u_f.  \label{ot4}
\end{align}
Therefore, the state estimates $(\hat w_f, \hat w_s, \hat u_s, \hat u_f)$ can be transformed into the state estimates  $(\hat\rho_{ f}, \hat v_{ f}, \hat \rho_{ s}, \hat v_{ s})$. The following conclusion is reached.
\begin{thm}
	Consider the two-lane traffic ARZ model in \eqref{rv1}-\eqref{rv4} with boundary conditions \eqref{bc1}-\eqref{bc4}, initial conditions $\rho_{ f}(x,0), v_{ f}(x,0), \rho_{ s}(x,0), v_{ s}(x,0) \in {L}^\infty\left([0,L]\right)$, state estimates $(\hat \rho_s(x,t), \hat \rho_f(x,t), \hat v_s(x,t), \hat v_f(x,t))$ are obtained from the collocated observer design \eqref{ot1}-\eqref{obl} for $(\hat w_s, \hat w_f, \hat v_s, \hat v_f)$ and the invertible transformation between them is given in \eqref{ot1}-\eqref{ot4}. The output injection gains $\{P(x)\}$ and $\{Q(x)\}$ are obtained in \eqref{PM}, \eqref{QN} by solving the kernels $\{M\}$ and $\{N\}$ from \eqref{mn1}-\eqref{mn2}. The finite-time convergence of estimation errors to zero equilibrium is reached in $t_o$ given by \eqref{to}.
\end{thm}
\begin{pf*}{Proof.}
	Lemma 3 with the existence of the backstepping transformation for the observer in \eqref{bkst1}, \eqref{bkst2} yields the convergence of estimation errors $(\check w_s,\check w_f,\check v_s,\check v_f)$ defined by \eqref{err1}-\eqref{errl} to zero for $t>t_o$. Given the transformation in \eqref{ot1}-\eqref{ot4}, the finite-time convergence to zero equilibrium is arrived for the estimation errors of the original system $(\tilde \rho_s(x,t), \tilde \rho_f(x,t), \tilde v_s(x,t), \tilde v_f(x,t))$. 
\end{pf*}

\subsection{Output feedback controller}
The output feedback controllers are constructed by employing the states estimates in the full-state feedback laws which yield the finite-time stability of the closed-loop system to zero equilibrium. Combining the collocated observer design \eqref{ob1}-\eqref{obl} and full-state feedback controllers \eqref{Us},\eqref{Uf}, we obtain the following output feedback controllers, 
\begin{align}
\notag U_s(t)=&\exp\left(-\frac{a^{vv}_{11}}{\mu_1}L\right) \int_{0}^{L} \frac{\gamma p_s^\star}{\rho_s^\star} K_{11}(L,\xi) \left(\hat\rho_{ s}(\xi, t)-\rho_{ s}^\star\right) 
\\ \notag  &\!\!\!\!\!+ \frac{\gamma p_f^\star}{\rho_f^\star}K_{12}(L,\xi) \left(\hat\rho_{ f}(\xi, t)-\rho_{ f}^\star\right)
\\ \notag &\!\!\!\!\! +\!\left[\!\!K_{11}(L,\xi)+L_{11}(L,\xi)\exp\left(\!\!\frac{a_{11}^{vv}}{\mu_1}\xi\!\right)\!\!\right]\!\! \left(\hat v_{ s}(\xi, t)- v_{ s}^\star\right) \\
&\!\!\!\!\!+\!\left[\!\!K_{12}(L,\xi)+L_{12}(L,\xi)\exp\left(\!\!\frac{a_{22}^{vv}}{\mu_2}\xi\!\right)\!\!\right]\!\! \left(\hat v_{ f}(\xi, t)-v_{ f}^\star\right)\!d\xi,\label{out1}\\
\notag U_f(t)=&\exp\left(-\frac{a^{vv}_{22}}{\mu_2}L\right)\int_{0}^{L} \frac{\gamma p_s^\star}{\rho_s^\star} K_{21}(L,\xi)  \left(\hat \rho_{ s}(\xi, t)-\rho_{ s}^\star\right) 
\\ \notag  &\!\!\!\!\!+\frac{\gamma p_f^\star}{\rho_f^\star}  K_{22}(L,\xi) \left(\hat \rho_{ f}(\xi, t)-\rho_{ f}^\star\right) 
\\ \notag  &\!\!\!\!\!+\!\left[\!\!K_{21}(L,\xi)+L_{21}(L,\xi)\exp\left(\!\!\frac{a_{11}^{vv}}{\mu_1}\xi\!\right)\!\!\right]\!\! \left(\hat v_{ s}(\xi, t)-v_{ s}^\star\right)\\
\notag	&\!\!\!\!\!+\!\left[\!\!K_{22}(L,\xi) + L_{22}(L,\xi)\exp\left(\!\!\frac{a_{22}^{vv}}{\mu_2}\xi\!\right)\!\!\right]\!\! \left(\hat v_{ f}(\xi, t)-v_{ f}^\star\right)\!d\xi,\label{out2}
\end{align}	  
where $(\hat\rho_{ f}, \hat \rho_{ s}, \hat v_{ s}, \hat v_{ f})$ are obtained from  $(\hat w_s, \hat w_f, \hat v_s, \hat v_f)$ using transformation in \eqref{ot1}-\eqref{ot4}.
\begin{thm}
	Consider the two-lane traffic ARZ model in \eqref{rv1}-\eqref{rv4} with boundary conditions \eqref{bc1}-\eqref{bc4}, initial conditions $\rho_{ f}(x,0), v_{ f}(x,0), \rho_{ s}(x,0), v_{ s}(x,0) \in {L}^\infty\left([0,L]\right)$ and the output feedback laws in \eqref{out1},\eqref{out2}, where the kernels $\{K\}$ and $\{L\}$ are obtained by solving \eqref{kn1}-\eqref{kn16} and output injection gains obtained by solving the kernels $\{M\}$ and $\{N\}$ in \eqref{mn1}-\eqref{mn16}. The steady states $(\rho_{ f}^\star, v_{ f}^\star, \rho_{ s}^\star, v_{ s}^\star)$ are finite-time stable and the convergence is reached in $t_{\rm out}$ defined as
	\begin{align}
		t_{\rm out} = t_o + t_f,
	\end{align}
	where $t_o$ is given in \eqref{to} and $t_f$ in \eqref{tf}.
\end{thm}
\begin{pf*}{Proof.}
	Theorem 3 yields that state estimates $(\hat\rho_{ f}, \hat \rho_{ s}, \hat v_{ s}, \hat v_{ f})$ converge to $(\rho_{ f}, \rho_{ s}, v_{ s}, v_{ f})$ after $t = t_o$. Applying Theorem 2, one has that $(\rho_{ f}, \rho_{ s}, v_{ s}, v_{ f})$ converge to  $(\rho_{ f}^\star, \rho_{ s}^\star, v_{ s}^\star, v_{ f}^\star)$ after $t = t_f$. Therefore, after $t =  t_o + t_f$, we have the convergence of state variables to steady states.
\end{pf*}
\section{Numerical Simulation}
To validate our control design including the full-state feedback controllers and the collocated boundary observer, we perform the numerical simulation for the two-lane ARZ model under two different secenrios of traffic congestion. For the first scenario, we consider the stop-and-go traffic appearing in the freeway segment of interest and therefore implement sinusoid initial conditions. For second scenario, we consider a single shock wave front for the initial state of traffic where the upstream vehicles are blocked by denser traffic downstream. This is a common phenomenon when slow moving vehicles block the road or changes of local road situations like hills and curves. Traffic bottleneck forms as a result. It follows the appearance of a moving shock wave consisting of high-density traffic downstream and relative low-density traffic upstream on the road. 

\begin{table}[t!]
	\caption{Model Parameter Table}
	\label{Table}
	\centering
	\begin{tabular}{|p{3.8cm}|p{1cm}|p{1cm}|p{1cm}|}
		\hline
		Name & Slow lane & Fast lane & Unit \\
		\hline
		Freeway segment length $L$ & $1000$ & $1000$ & $\unit{m}$ \\
		\hline
		Pressure exponent $\gamma$  &$0.8$ & $0.8$ & $\unit{1}$ \\
		\hline
		Maximum density $\rho_i^m$ & $240$ & $150$ & $\unit{veh}/\unit{km}$ \\
		\hline
		Maximum velocity $v_m$ & $144$ & $144$ & $\unit{km}/\unit{h}$ \\
		\hline
		Steady state densities $\rho_i^\star$  & $180$ & $80$ & $\unit{veh}/\unit{km}$ \\
		\hline
		Steady state velocity $v_i^\star$ & $32$ & $40$ & $\unit{km}/\unit{h}$ \\
		\hline
		Relaxation time $T_i^e$ & $200$ & $100$ & $\unit{s}$ \\ 
		\hline 
		Driver's lane preference $T_i$ & $50$ & $25$ & $\unit{s}$ \\ 
		\hline 
	\end{tabular}
\end{table}

The control design presented in the previous sections are tested and illustrated for both secenrios. The model parameters used in the numerical simulation are given in Table~\ref{Table}. Both the fast-lane and slow-lane are considered in the congested regime where the vehicles on the road are relatively dense so that the velocity disturbances propagate from the leading vehicle to the following vehicle. Steady states density $\rho_{ i}^\star$ are chosen given the maximum density $\rho_i^m$ and Maximum velocity $v_m$ so that the traffic of both lanes are lightly congested. We consider the situation that in general drivers prefer the slow lane rather than the fast lane. $T_f$ is smaller than $T_s$ since drivers prefer remaining in the slow lane rather than changing to the fast lane. Therefore, higher density traffic appears in the slow lane and it can contain higher traffic flow. Steady state velocity $v_{ i}^\star$ are obtained based on this parameter choice.

In the following figures, the evolution of the state variables are illustrated with surface plots. The initial conditions of the states are highlighted with color blue and the outlet boundary control inputs are highlighted with color red.

\begin{figure}[t!]
	\includegraphics[width=10cm]{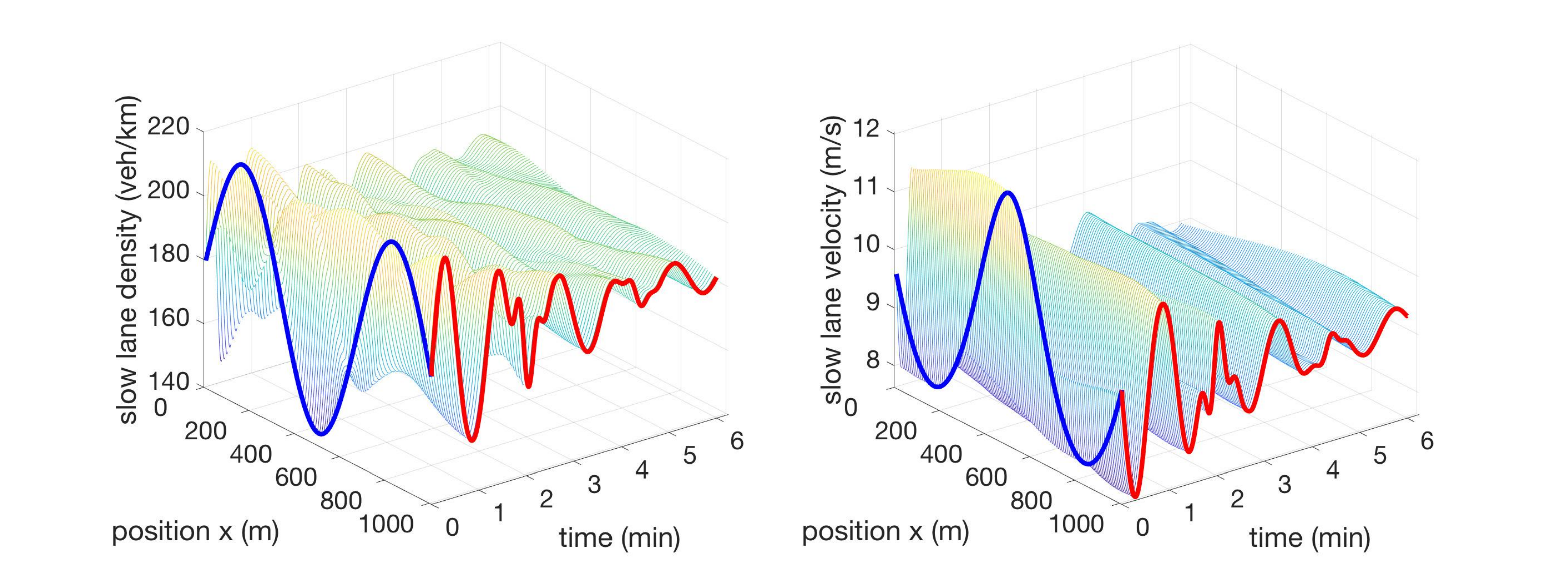}
	\caption{Scenario 1: density and velocity of slow lane traffic of open-loop system with sinusoid initial conditions.}
	\label{fig:open_s}
	\includegraphics[width=10cm]{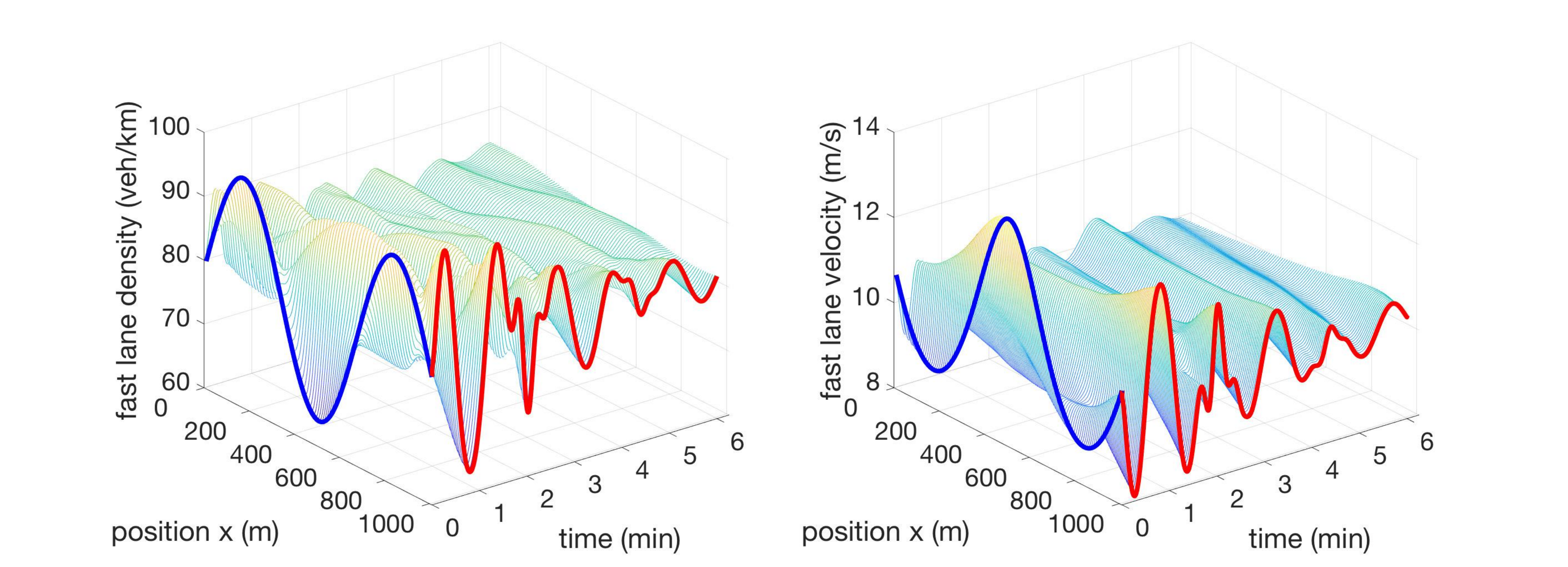}
	\caption{Scenario 1: density and velocity of fast lane traffic of open-loop system with sinusoid initial conditions.}
	\label{fig:open_f}
\end{figure}
\begin{figure}[t!]
	\includegraphics[width=10cm]{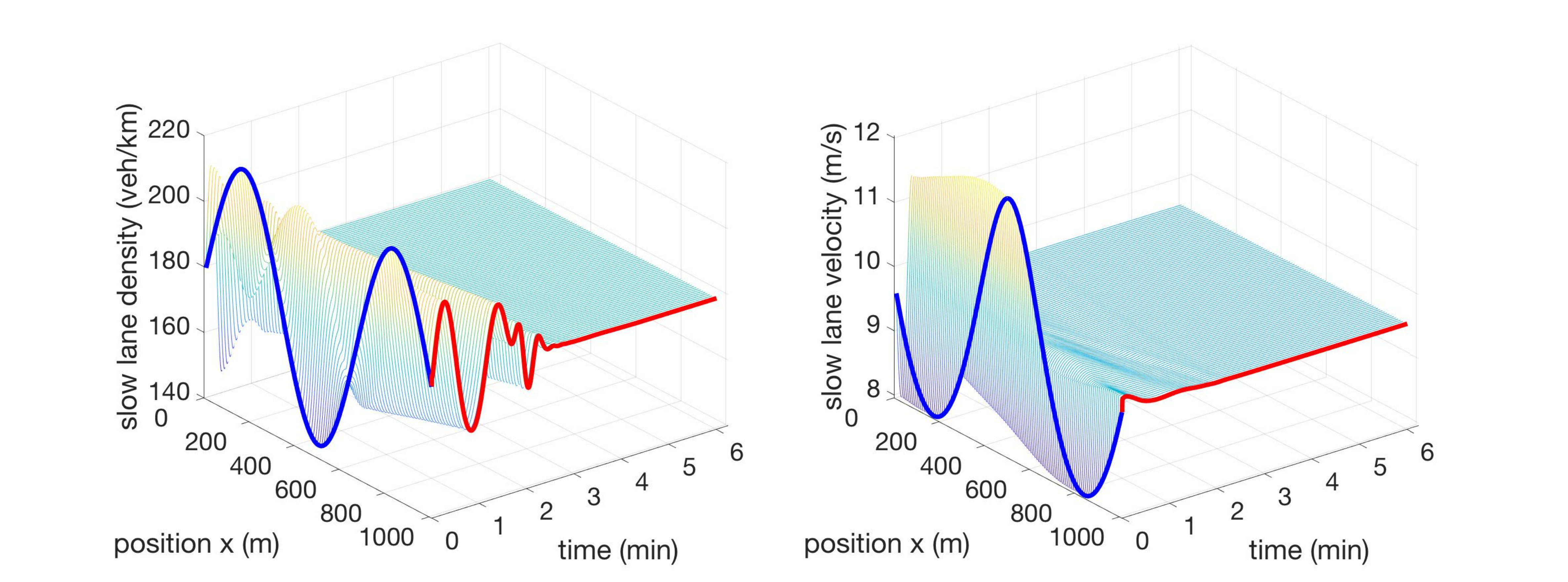}
	\caption{Scenario 1: density and velocity of slow lane traffic of closed-loop system with full-state feedback controllers.}
	\label{fig:full_s}
	\includegraphics[width=10cm]{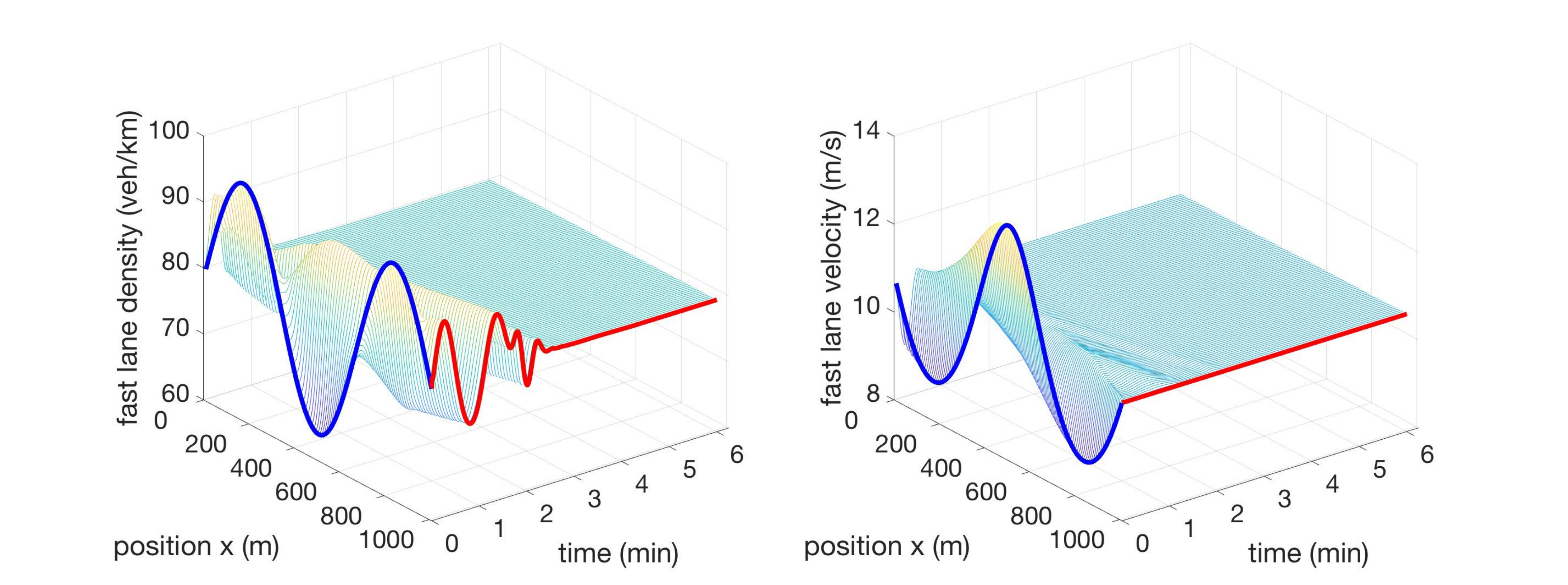}
	\caption{Scenario 1: density and velocity of fast lane traffic of closed-loop system with full-state feedback controllers.}
	\label{fig:full_f}	
\end{figure}

\subsection{Scenario 1: stop-and-go traffic}
In this scenario, we consider the initial traffic states are oscillated around the equilibrium states and thus we implement sinusoid initial conditions for traffic density and velocity. The constant incoming flow and outgoing flow are considered for the open-loop simulation as shown in Fig.~\ref{fig:open_s} and Fig.~\ref{fig:open_f}. For the steady state velocity, it takes around $100\;\rm s$ for both fast-lane and slow-lane vehicles to leave the considered freeway segment. But the oscillations sustain for more than $6\;\rm \min$. 

The full state feedback stabilization results are shown in Fig.~\ref{fig:full_s} and Fig.~\ref{fig:full_f}. The finite-time convergence of the density and velocity states to the steady states is achieved in $t_f = 260 \;\rm s$.

\begin{figure}[t!]
	\includegraphics[width=10cm]{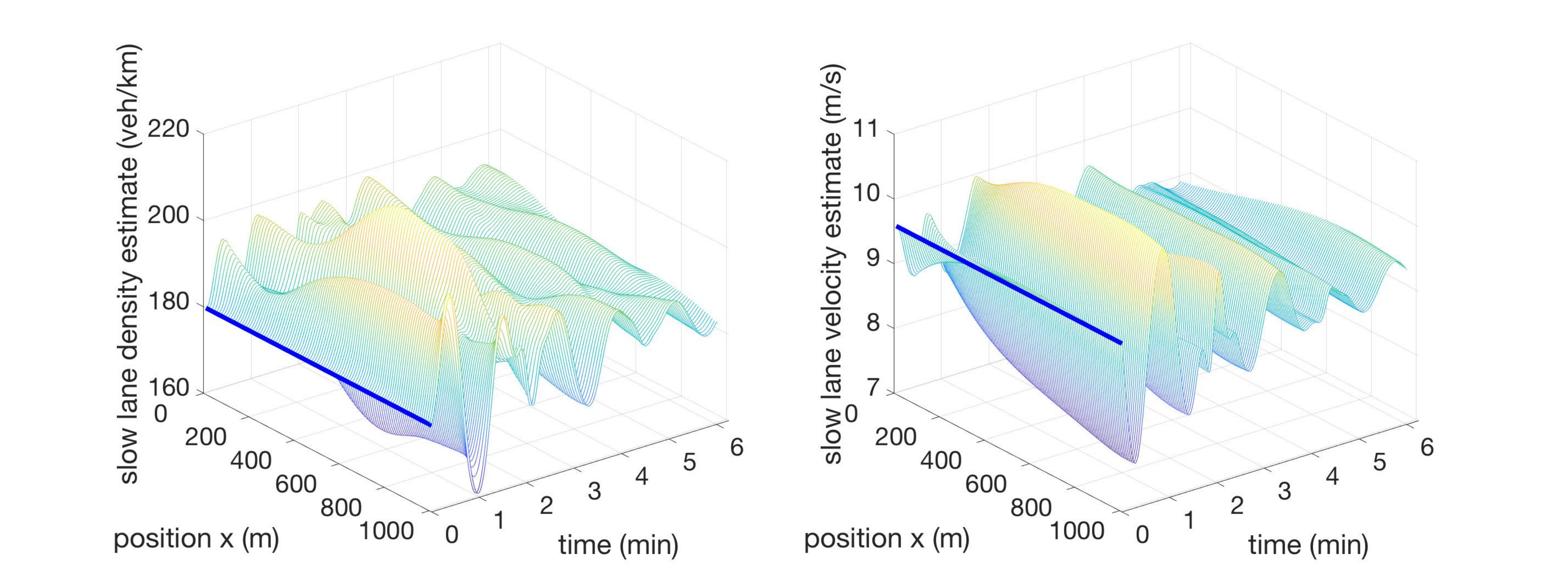}
	\caption{Scenario 1: density and velocity estimates of slow lane traffic of open-loop system with sinusoid initial conditions.}
	\label{fig:open_se}
	\includegraphics[width=10cm]{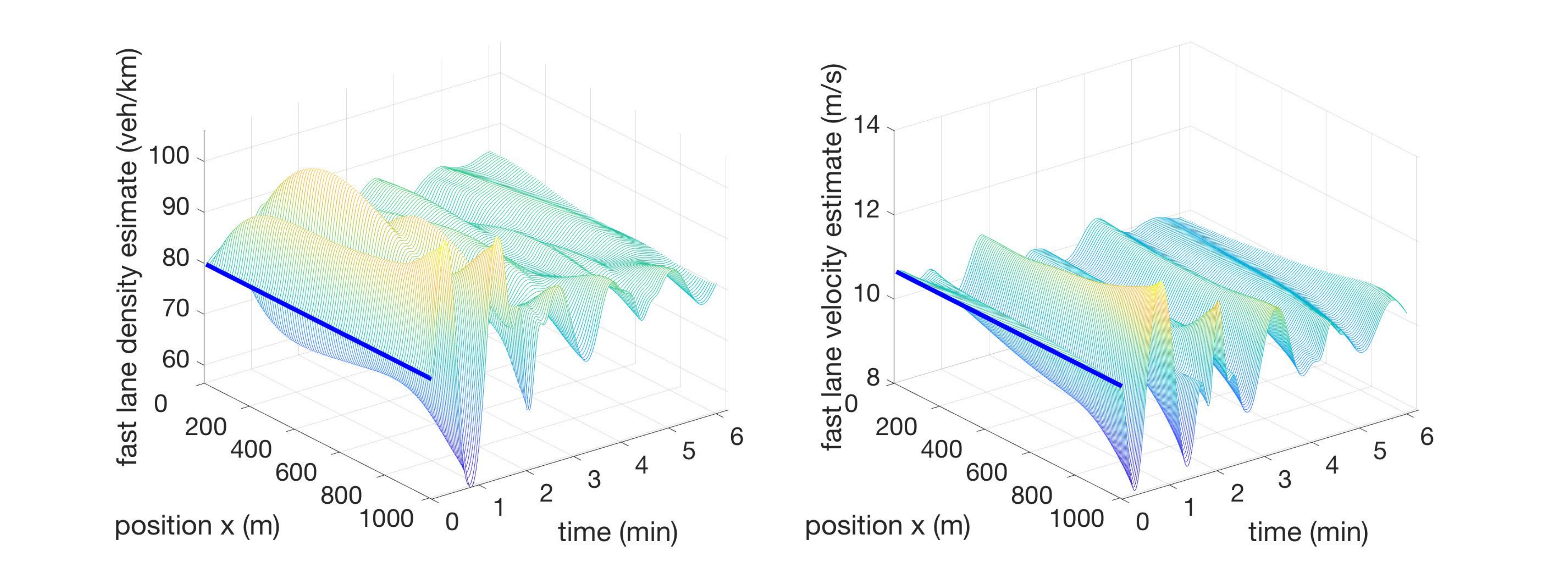}
	\caption{Scenario 1: density and velocity estimates of fast lane traffic of open-loop system with sinusoid initial conditions.}
	\label{fig:open_fe}
\end{figure}

\begin{figure}[t!]
	\includegraphics[width=9.7cm]{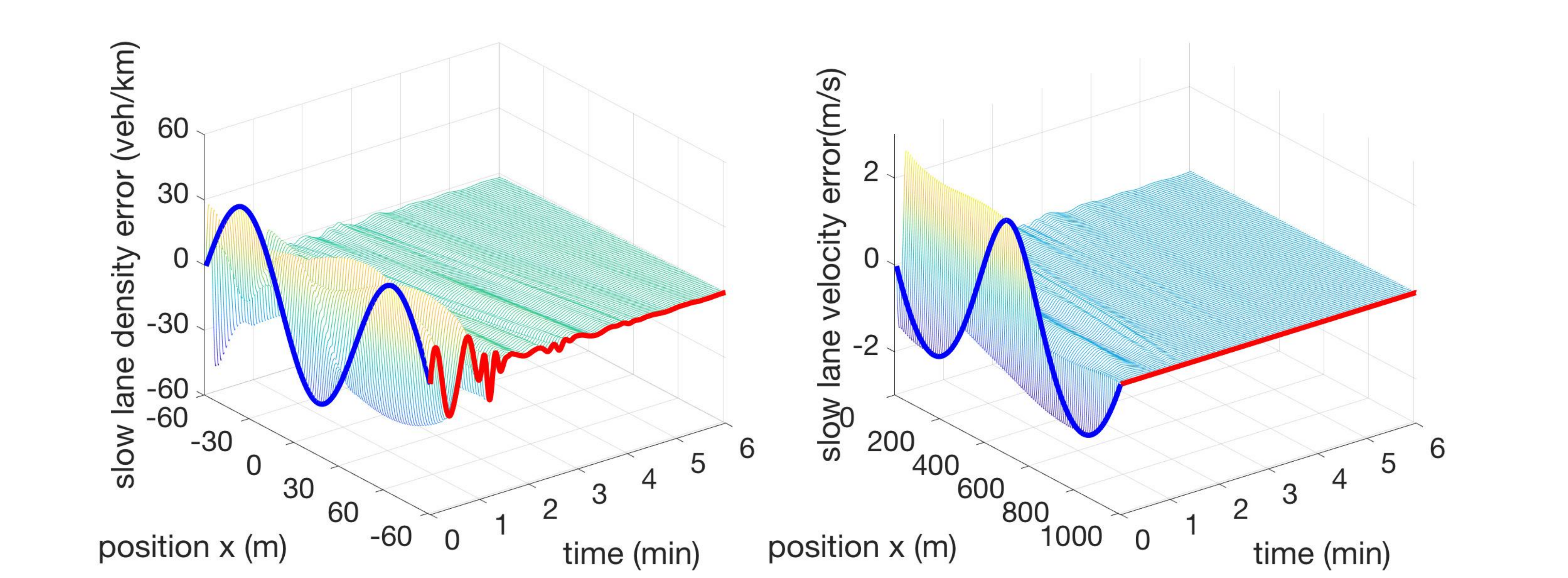}
	\caption{Scenario 1: density and velocity estimation errors of slow lane traffic of closed-loop system with full-state feedback controllers.}
	\label{fig:err_f} 
	\includegraphics[width=9.7cm]{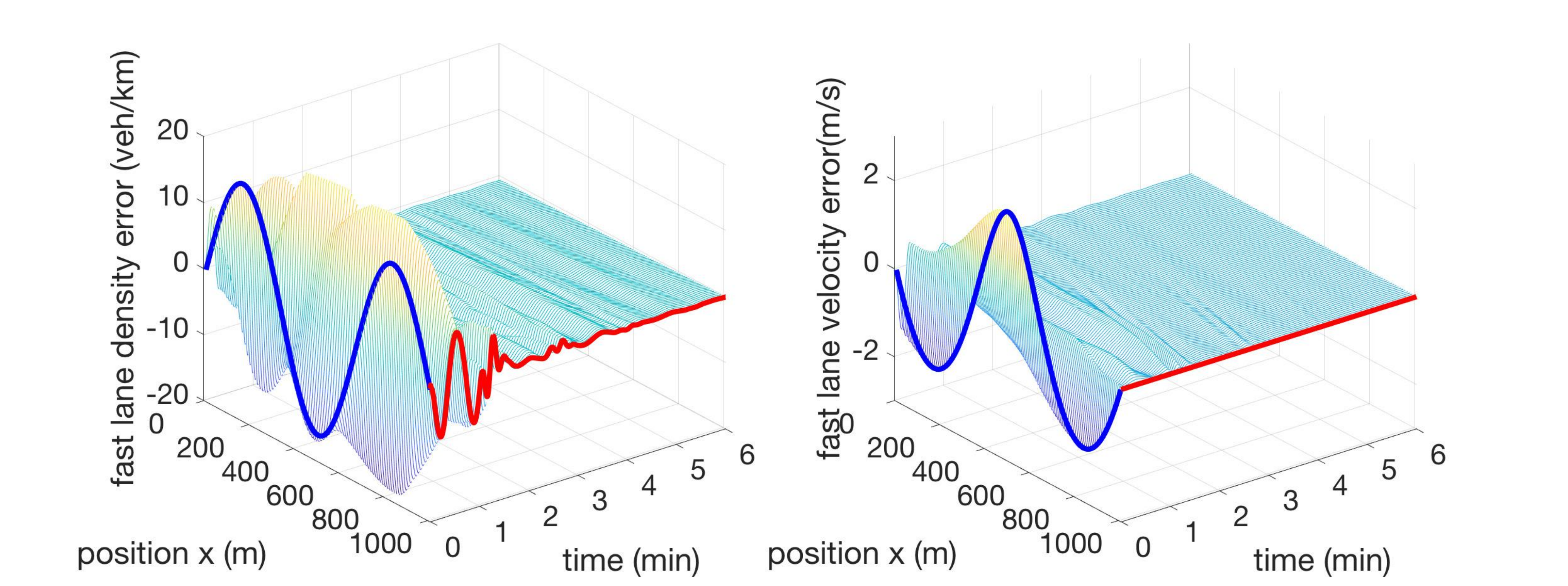}
	\caption{Scenario 1: density and velocity estimation errors of fast lane traffic of open-loop system with full-state feedback controllers.}
	\label{fig:err_s}
\end{figure}

The simulation results of the collocated observer design is shown in Fig.~\ref{fig:open_se} and Fig.~\ref{fig:open_fe}. The estimation errors are plotted in Fig.~\ref{fig:err_f} and Fig.~\ref{fig:err_s}. Here we choose open-loop system to validate observer design since the oscillations are not damped out to zero by control design and therefore estimation result could be illustrated better in this case. Without knowledge of the initial state of the system, we implement uniform steady state value for initial conditions, highlighted with blue lines. From Fig.~\ref{fig:err_f} and Fig.~\ref{fig:err_s}, we can see that after $t_o = 310 \;\rm s$, the estimation errors converge to zero, indicating that the state estimates in Fig.~\ref{fig:open_se} and Fig.~\ref{fig:open_fe} converge to the open-loop simulation of the states in Fig.~\ref{fig:open_s} and Fig.~\ref{fig:open_f}.

\begin{figure}[t!]
	\centering
	\includegraphics[width=10cm]{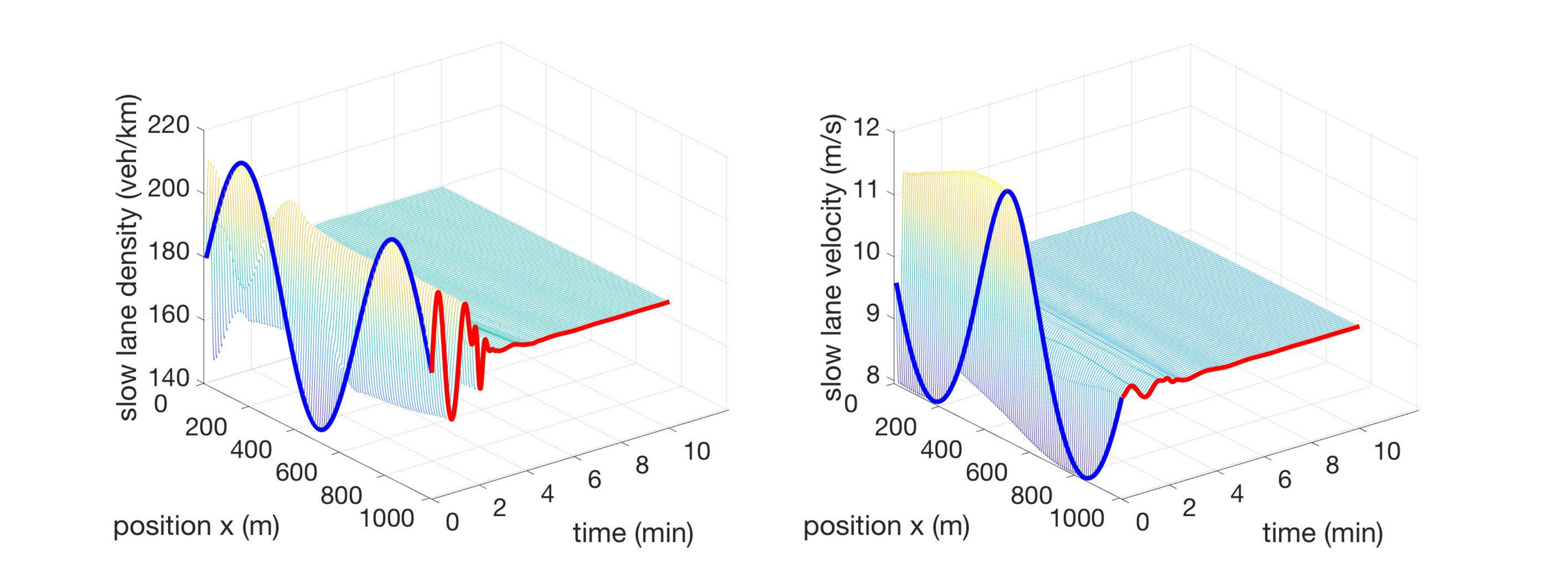}
	\caption{Scenario 1: density and velocity of slow lane closed-loop system with output feedback controllers.}
	\label{fig:output_s}
	\includegraphics[width=10cm]{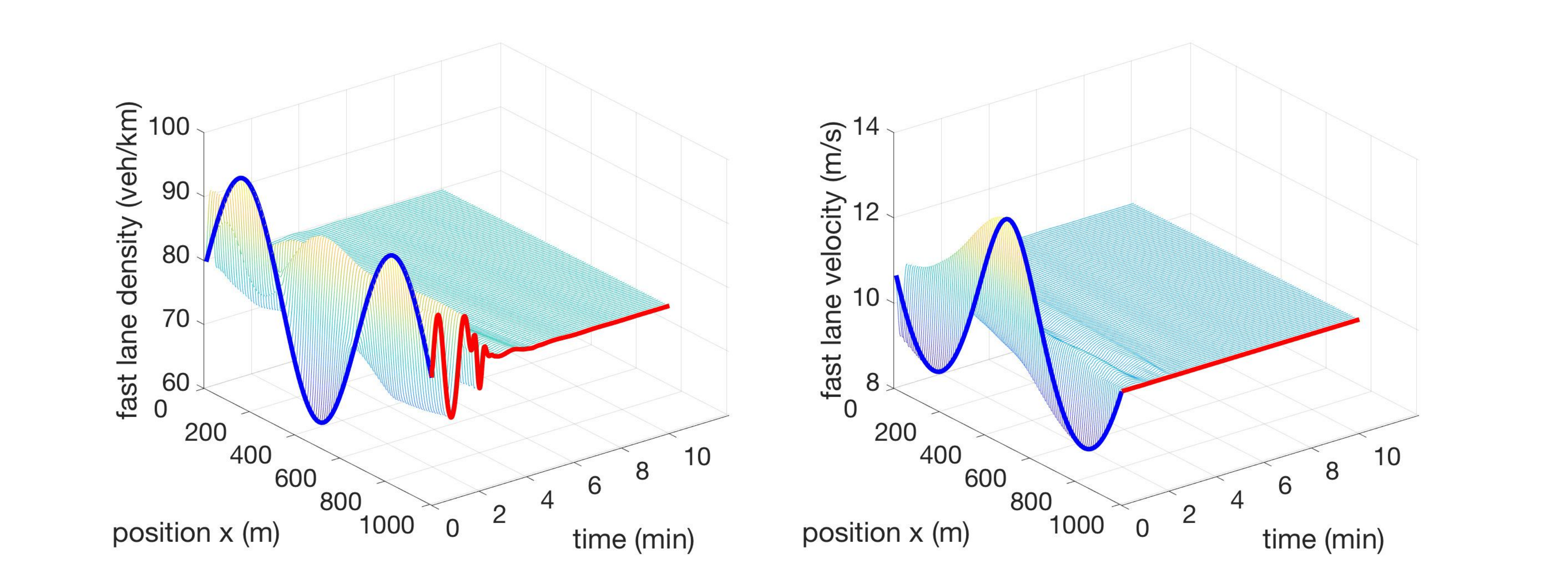}
	\caption{Scenario 1: density and velocity of fast lane closed-loop system with output feedback controllers.}
	\label{fig:output_f}	
\end{figure}

Combining the observer design and the full-state feedback controllers, we derive the output feedback controllers and then simulate the closed-loop system in Fig.~\ref{fig:output_s} and Fig.~\ref{fig:output_f}. The finite convergence time of the closed-loop with output feedback controllers are $t = t_o + t_f = 570 \rm s$. It is shown in the figures that the states converge to the steady state values before $10\;\rm \min$.

\begin{figure}[t!]
	\centering
	\includegraphics[width=10cm]{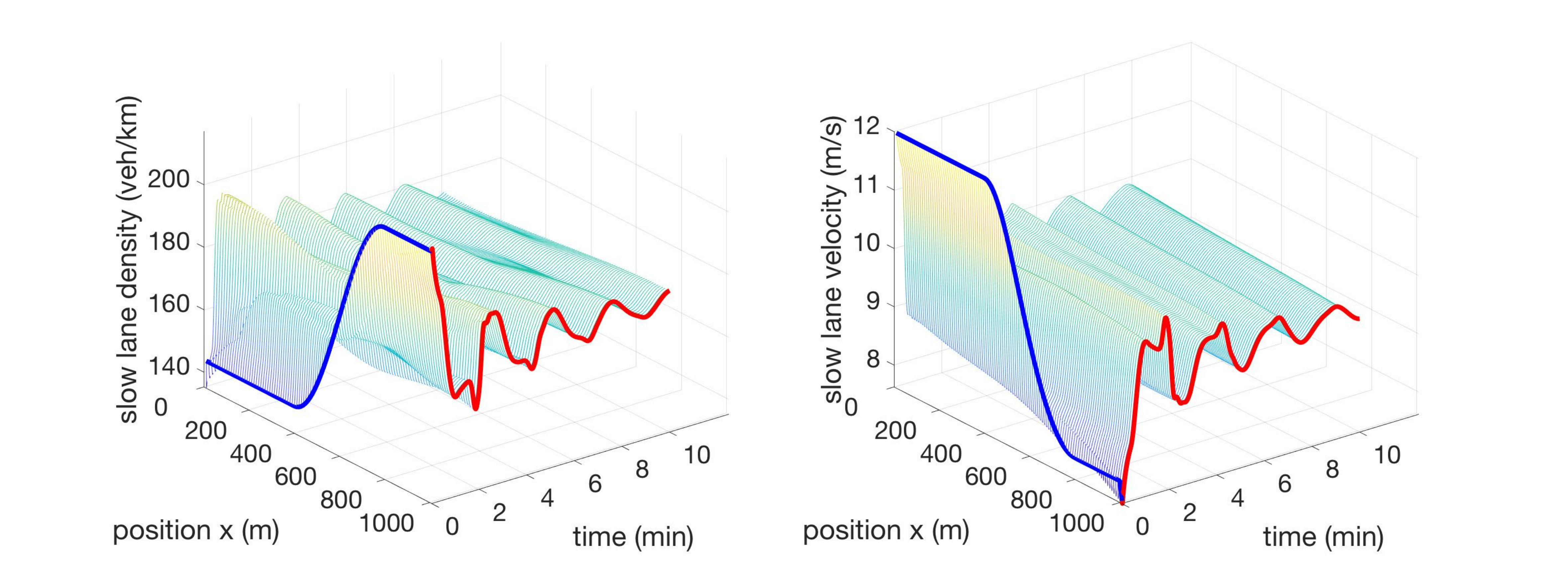}
	\caption{Scenario 2: density and velocity of slow lane traffic of open-loop system with shockwave initial conditions.}
	\label{fig:SSopens}
	\includegraphics[width=10cm]{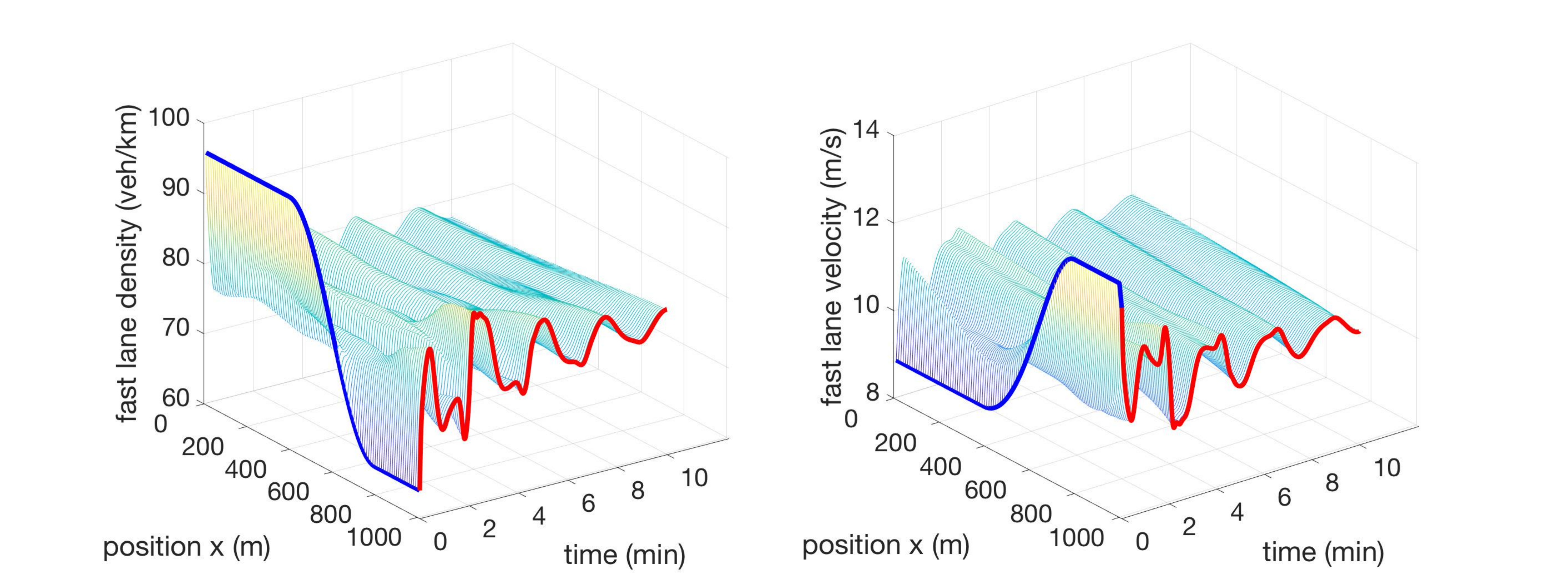}
	\caption{Scenario 2: density and velocity of fast lane traffic of open-loop system with shockwave initial conditions.}
	\label{fig:SSopenf}
\end{figure}	

\subsection{Scenario 2: traffic bottleneck}
Consider in Scenario 2 that there are some local changes of road situations like uphill and downhill gradients, curves downstream of the freeway segment. Therefore, traffic bottleneck forms from the downstream. We implement a shockwave front shape of the initial conditions for the slow-lane where traffic densities close to the outlet of the segment are denser and light densities traffic is blocked at the upstream. As a result, the traffic velocity is faster near the inlet while the velocity become slower near the outlet. On the other hand, we consider for the fast-lane that there is traffic flow of high density entering from the inlet. In general, drivers prefer the slow lane and the relative light traffic flow close to the inlet will trigger the lane-changing from the fast lane to slow lane close to the inlet. This worsens the traffic congestion on the slow lane since the traffic bottleneck appears in the downstream of slow lane.  

\begin{figure}[t!]
	\centering
	\includegraphics[width=10cm]{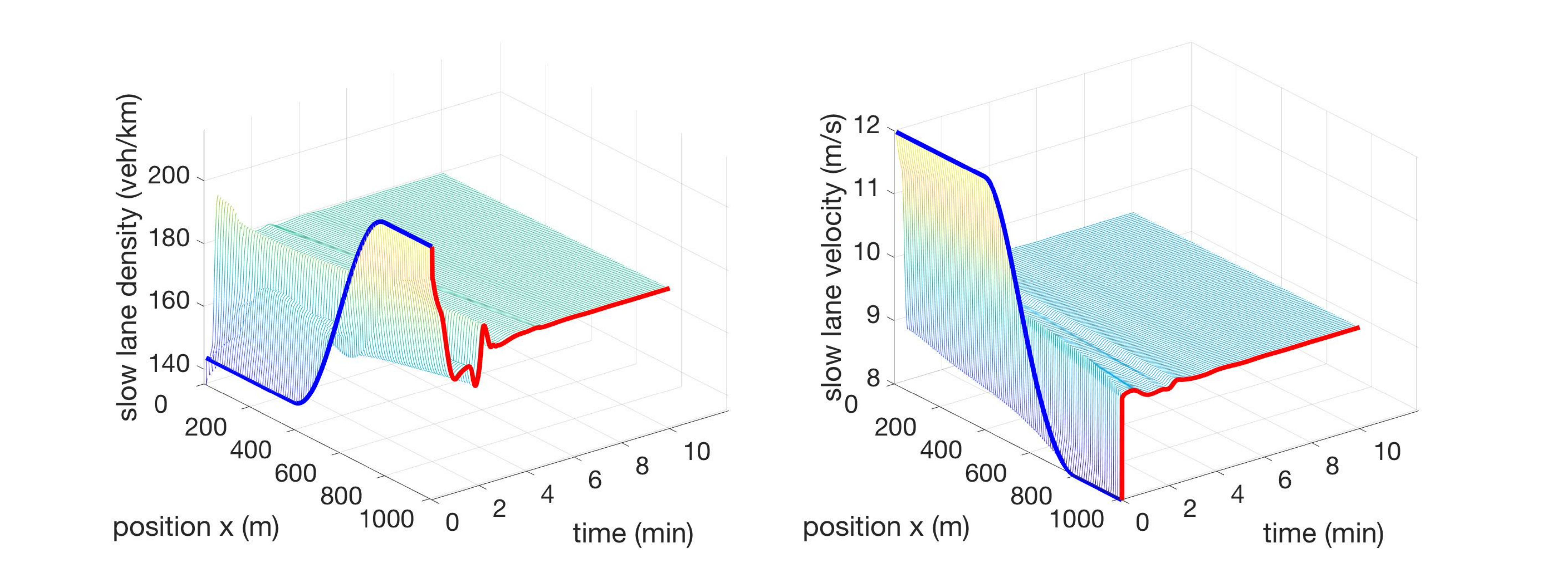}
	\caption{Scenario 2: density and velocity of slow lane closed-loop system with output feedback controllers.}
	\label{fig:SS_outputs}
	\includegraphics[width=10cm]{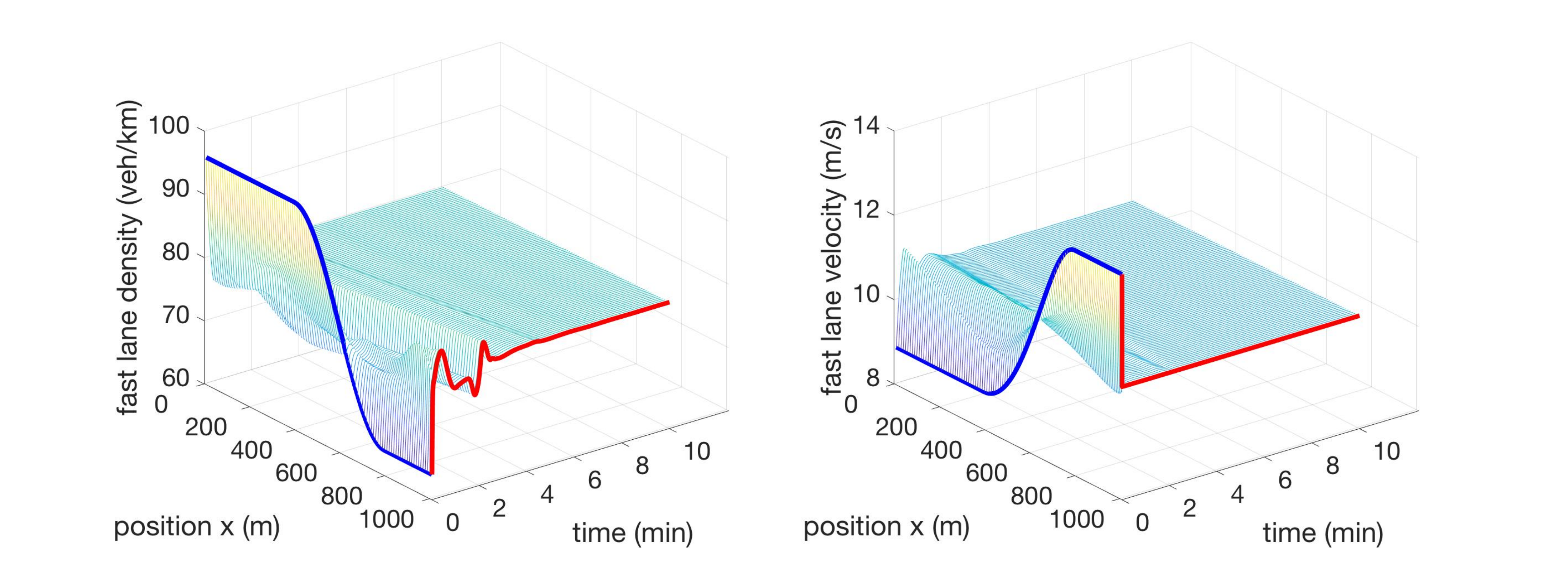}
	\caption{Scenario 2: density and velocity of fast lane closed-loop system with output feedback controllers.}
	\label{fig:SS_outputf}
\end{figure}

In the open-loop simulation shown by Fig.~\ref{fig:SSopens} and Fig.~\ref{fig:SSopenf}, soft shock wave initial traffic states result in the stop-and-go traffic on the freeway segment. Here we omit the state estimation results by collocated observer which has been demonstrated in Scenario 1. The simulation result of the output feedback control applied to Scenario 2 is given with Fig.~\ref{fig:SS_outputs} and Fig.~\ref{fig:SS_outputf}. We can see that the oscillated traffic congestion of $1 \rm km$ in Fig.~\ref{fig:SSopens} and Fig.~\ref{fig:SSopenf} is damped out in a fast manner for around $4\;\rm \min$.

With the numerical simulation of the two-lane ARZ model with lane changing in two different secenrios, we demonstrate that the full-state feedback controllers, the collocated observer and the output feedback controllers achieve the finite-time convergence of the state variations from the steady states and estimation errors to zero.

\section{Conclusion}
This paper solves the output feedback stabilization of a two-lane traffic congestion problem with lane-changing. Using coordinate transformation and backstepping method, the linearized first-order coupled $4 \times 4$ hyperbolic PDE system is transformed into a cascade target system. The finite-time convergence to the steady states is achieved with two VSLs control inputs actuating velocities at the outlet. By taking measurement of density variations at the outlet, a collocated observer design is proposed for state estimation, which is theoretically novel and practically sound. This result paves the way for applying PDE backstepping techniques for multi-lane traffic with inter-lane activities. There are concerns on modeling the lane-changing as density exchanging source terms. However, the control design and the methodology proposed in this paper should not be limited by possible modifications on these terms.

\end{document}